\documentclass[a4paper,11pt]{article}
\usepackage{amsmath, amsfonts, amssymb}
\usepackage[english, russian]{babel}
\numberwithin{equation}{section}
\newtheorem{lemma}{Лемма}
\newtheorem{theorem}{Теорема}

\newtheorem{corollary}{Следствие}

\usepackage{graphicx}

\allowdisplaybreaks  

\begin{document}

\begin{center}
Double- and simple-layer potentials for generalized singular
elliptic equations and their applications to the solving the
Dirichlet problem
\end{center}

\begin{center}
\textbf{Ergashev T.G.}
\end{center}

Institute of Mathematics 81, Mirzo Ulugbek street, Tashkent,
100170 Uzbekistan

Institute of Irrigation and Agricultural Mechanization Engineers
39, Kari-Niyazi Street, Tashkent 700100 Uzbekistan

E-mail: ergashev.tukhtasin@gmail.com
\\

\textbf{Abstract.} Potentials play an important role in solving
boundary value problems for elliptic equations. In the middle of
the last century, a potential theory was constructed for a
two-dimensional elliptic equation with one singular coefficient.
In the study of potentials, the properties of the fundamental
solutions of the given equation are essentially and fruitfully
used. At the present time, fundamental solutions of a
multidimensional elliptic equation with several singular
coefficients are already known. In this paper, we investigate the
double- and simple-layer potentials for this kind of elliptic
equations. Results from potential theory allow us to represent the
solution of the boundary value problems in integral equation form.
By using a decomposition formula and other identities for the
Lauricella's hypergeometric function in many variables, we prove
limiting theorems and derive integral equations concerning a
densities of the double- and simple-layer potentials. The obtained
results  are applied to find an explicit solution of the Dirichlet
problem for the generalized singular elliptic equation in the some
part of the multidimensional ball.



\section{Introduction}

Potential theory has played a paramount role in both analysis and
computation for boundary value problems for elliptic equations.
Numerous applications can be found in solid mechanics, fluid
mechanics, elastic dynamics, electro-magnetics, and acoustics.
Results from potential theory allow us to represent boundary value
problems in integral equation form. For problems with known
Green's functions, an integral equation formulation leads to
powerful numerical approximation schemes.

The double- and simple-layer potentials play an important role in
solving boundary value problems for elliptic equations. For
example, the representation of the solution of the Dirichlet
problem  for the Laplace equation is sought as a double-layer
potential with unknown density and an application of certain
property leads to a Fredholm equation of the second kind for
determining the density function (see \cite{Gun} and \cite{Mir}).

Interest in the potential theory for the singular elliptic
equation has increased significantly after Gellerstedt's papers
\cite{{Gel1},{Gel2}}. In works \cite{Frankl} and \cite{Pul}, the
potential theory was exposed for the following  simplest
degenerating elliptic equation
\begin{equation}
\label{simplest} {{\frac{\partial^2u}{\partial
x_1^2}}}+{{\frac{\partial^2u}{\partial x_2^2}}} +
{\frac{{2\alpha_1}} {{x_{1}}} }{{\frac{\partial u}{\partial x_1}}}
= 0, \,\,0 < 2\alpha_1 < 1
\end{equation}
in the domain, which is bounded in the half-plane $x_1>0$. An
exposition of the results on the potential theory for the
two-dimensional singular elliptic equation (\ref{simplest})
together with references to the original literature are to be
found in the monograph by Smirnov \cite{Sm}, which is the standard
work on the subject. This work also contains an extensive
bibliography of all relevant papers up to 1966; the list of
references given in the present work is largely supplementary to
Smirnov's bibliography. Various interesting problems associated
with the equation (\ref{simplest}) were studied by many authors
(see \cite{{A1}, {A2},{AY},{G5},{G4},{G7},{G9a},{He},{Lo},{NL}}).

In his work \cite{Has} Hasanov found fundamental solutions of the
generalized bi-axially symmetric Helmholtz equation
\begin{equation}
\label{2Helmholtz} {{\frac{\partial^2u}{\partial
x_1^2}}}+{{\frac{\partial^2u}{\partial x_2^2}}} +
{\frac{{2\alpha_1}} {{x_{1}}} }{{\frac{\partial u}{\partial x_1}}}
+ {\frac{{2\alpha_2}} {{x_{2}}} }{{\frac{\partial u}{\partial
x_2}}}+\lambda u = 0,
\end{equation}
in the domain, which is bounded in the quarter-plane $x_1>0$,
$x_2>0$, where $\alpha_1,\,\alpha_2$ and $\lambda$ are real
numbers ($0 < 2\alpha_1,\,2\alpha_2 < 1$).   When $\lambda=0$,
this equation is known as the equation of the generalized
bi-axially symmetric potential theory whose name is due to
Weinstein who first considered fractional dimensional space in
potential theory \cite{{W48},{W53}}. Using Hasanov's results
authors of the works \cite{{Berd},{Tomsk},{Ufa},{SHC}} constructed
the double-layer potential theory for the two-dimensional equation
(\ref{2Helmholtz}) in the case when $\lambda=0$.

Relatively few papers have been devoted to the potential theory
and boundary value problems for the singular elliptic equations
when the dimension exceeds two \cite{{Hu},{Mav},{Muh}}.

In the present work we shall give the potential theory for the
following generalized singular elliptic equation
\begin{equation}
\label{eq1} H_{\alpha} ^{(m,n)} (u) \equiv {\sum\limits_{i =
1}^{m} {u_{x_{i} x_{i}}} }  +
\sum\limits_{k=1}^n{\frac{{2\alpha_k}} {{x_{k}}} }u_{x_{k}}  = 0,
\end{equation}
where $m$ is a dimensional of the Euclidean space $R_m$, $n$ is a
number of the singular coefficients of elliptic equation $(m \ge
2$,\,\,$0 < n \le m);$ $\alpha=(\alpha_1,...,\alpha_n)$ and $
\alpha _{k} $ are real numbers with $0 < 2\alpha _{k} < 1$, $k\in
K$\, $(K=\{1,...,n\})$, and apply this theory to the finding a
regular solution of the Dirichlet problem for  equation
(\ref{eq1}) in the domain, which is  bounded in the subset of the
Euclidean space
\begin{equation*}
R_m^{n+}=\left\{x=\left(x_1,...,x_m\right) \in R_m:
x_1>0,...,x_n>0\right\}.
\end{equation*}

By a regular solution of the equation (\ref{eq1}) is meant a
function that has continuous derivatives up to the second order
(inclusive) in some domain and satisfies the equation (\ref{eq1})
at all points of this domain.

Naturally, in solving the problem posed for the equation
(\ref{eq1}), an important role is played some fundamental solution
of this equation. Fundamental solutions of equation (\ref{eq1})
were constructed recently \cite{A28}. In fact, the fundamental
solutions of generalized singular elliptic equation can be
expressed in terms of Lauricella's hypergeometric function in many
variables  $F_A^{(n)}(a, b_1,...,b_n; c_1,...,c_n; y_1,...,y_n)$
(see, for details, \cite{A30}) and it becomes clear that the
number of variables of a hypergeometric function is equal to the
number of singular coefficients of the equation (\ref{eq1}).
Therefore, in order to facilitate the constructing process  of the
potential theory for equation (\ref{eq1}), we preliminary study
some necessary properties of this Lauricella's hypergeometric
function in many variables.

\section{Preliminaries}

With a view to introducing formally the Gaussian hypergeometric
function and its generalization (that is, Lauricella
hypergeometric function in several variables), we recall here some
definitions and identities involving Pochhammer's symbol
$(\lambda)_p$, and the Gamma function $\Gamma(z)$ ($\lambda$ and
$z$ are complex numbers) defined by
\begin{equation*} \label{gamma}  \Gamma(z)=
\left\{ {{\begin{array}{*{20}c}
 {\displaystyle\int_{0}^{\infty}t^{z-1}e^{-t}dt,\,\, Re(z)>0,} \hfill \\
 \\
 {\displaystyle\frac{\Gamma(z+1)}{z}, \,\,\,\,\,\,\,\,\,\,\,\,\,\,\,\,Re(z)<0}; \, z\neq -1,-2,-3,...\,. \hfill \\
\end{array}}} \right.
\end{equation*}

Throughout this work we shall find it convenient to employ the
Pochhammer symbol $(\lambda)_p$ defined by
\begin{equation*} \label{pochh}  (\lambda)_p=
\left\{ {{\begin{array}{*{20}c}
 {1,\,\,\,\,\,\,\,\,\,\,\,\,\,\,\,\,\,\,\,\,\,\,\,\,\,\,\,\,\,\,\,\,\,\,\,\,\,\,\,\,\,\,\,\,\,\,\,\,\,\,\,\,\,\,\,  \textrm{if}\,\,p=0,} \hfill \\
 {\lambda(\lambda+1)...(\lambda+p-1), \,\,\textrm{if}\,\,\, p \in \mathbb{N}}, \hfill \\
\end{array}}} \right.
\end{equation*}
where $\mathbb{N}:=\{1,2,3,...\}.$ Since $(1)_p=p!$, $(\lambda)_p$
may be looked upon as a generalization of the elementary
factorial; hence the symbol $(\lambda)_p$ is also referred to as
the factorial function.

The Gaussian hypergeometric function is defined in the form (see,
for details, \cite{{E},{SK}})
\begin{equation} \label{hyperg} F\left( {a,b;c;z} \right) \equiv F{\left[
{{\begin{array}{*{20}c}
 {a,b;} \hfill \\
 {c;} \hfill \\
\end{array}} z} \right]}: = {\sum\limits_{m = 0}^{\infty}  {{\frac{{\left(
{a} \right)_{m} \left( {b} \right)_{m}}} {{\left( {c} \right)_{m}
m!}}}z^{m},\,c \ne 0,-1,-2,...\,,}}
\end{equation}
where $a,b,c$ and $z$  are complex numbers.

It is easily seen that the hypergeometric function $F\left(
{a,b;c;z} \right) $ in (\ref{hyperg}) converges absolutely within
the circle, that is, when $|z|<1$, provided that the denominator
parameter $c$ is neither zero nor a negative integer. Further
tests readily show that the hypergeometric function in
(\ref{hyperg}), when $|z|=1$ (that is, on the unit circle), is
absolutely convergent if $Re(c-a-b)>0$.

The linear transformation of the hypergeometric function, known as
Euler's transformation, may be recalled here as follows  \cite{E}:
\begin{equation}
\label{trans} F\left( {a,b;c;z} \right) = \left( {1 - z} \right)^{
- b}F\left({c-a,b;c;{\frac{{z}}{{z - 1}}}}\right).
\end{equation}

The Lauricella's hypergeometric function in $n \in \mathbb{N}$
variables $y:=\left(y_1,...,y_n\right)$ has a form \cite{A30}
\\

$
 F_{A}^{(n)} {\left[
{{\begin{array}{*{20}c}
 {a,b_{1} ,...,b_{n} ;} \hfill \\
 {c_{1} ,...,c_{n} ;} \hfill \\
\end{array}} y}  \right]}\equiv F_{A}^{(n)} {\left[
{{\begin{array}{*{20}c}
 {a,b_{1} ,...,b_{n} ;} \hfill \\
 {c_{1} ,...,c_{n} ;} \hfill \\
\end{array}} y_1,...,y_n}  \right]}
$
\begin{equation}
\label{eq7}
 = {\sum\limits_{m_{1} ,...m_{n} = 0}^{\infty}  {{\frac{{\left( {a}
\right)_{m_{1} + ... + m_{n}}  \left( {b_{1}}  \right)_{m_{1}}
...\left( {b_{n}}  \right)_{m_{n}}} } {{\left( {c_{1}}
\right)_{m_{1}}  ...\left( {c_{n}}  \right)_{m_{n}}} }
}{\frac{{y_{1}^{m_{1}}} } {{m_{1} !}}}...{\frac{{y_{n}^{m_{n}}} }
{{m_{n} !}}}}},
\end{equation}
where $a, b_k$, $c_k$ and $y_k$ are complex numbers and $c_{k} \ne
0,-1,-2,...,\,k \in K$.

We give the elementary relations for $F_{A}^{(n)}$ necessary in
this study:
\\

$ {\displaystyle\frac{{\partial}} {{\partial y_{k}}}}F_{A}^{\left(
{n} \right)} {\left[ {{\begin{array}{*{20}c}
 {a,b_{1},...,b_{n};} \hfill \\
 {c_{1} ,...,c_{n};} \hfill \\
\end{array}} y}\right]}
$
\begin{equation}
\label{eq24001}=\frac{ab_k}{c_k}F_{A}^{\left( {n} \right)} {\left[
{{\begin{array}{*{20}c}
 {a+1,b_{1},...,b_{k-1},b_k+1,b_{k+1},...,b_{n};} \hfill \\
 {c_{1} ,...,c_{k-1},c_k+1,c_{k+1},...,c_{n};} \hfill \\
\end{array}} y}\right]},
\end{equation}
\\

$
 \displaystyle\sum\limits_{k = 1}^{n}
\displaystyle\frac{b_k}{c_k}y_kF_{A}^{\left( {n} \right)} {\left[
{{\begin{array}{*{20}c}
 {a+1,b_{1},...,b_{k-1},b_k+1,b_{k+1},...,b_{n};} \hfill \\
 {c_{1} ,...,c_{k-1},c_k+1,c_{k+1},...,c_{n};} \hfill \\
\end{array}} y}\right]}
$
\begin{equation}
\label{eq25}=F_{A}^{\left( {n} \right)} {\left[
{{\begin{array}{*{20}c}
 {a+1,b_{1},...,b_{n};} \hfill \\
 {c_{1} ,...,c_{n};} \hfill \\
\end{array}} y}\right]}-F_{A}^{\left( {n} \right)} {\left[ {{\begin{array}{*{20}c}
 {a,b_{1},...,b_{n};} \hfill \\
 {c_{1} ,...,c_{n};} \hfill \\
\end{array}} y}\right]},
\end{equation}
\\

$ \displaystyle\frac{ab_k}{\left(c_k-1\right)c_k}y_kF_{A}^{\left(
{n} \right)} {\left[ {{\begin{array}{*{20}c}
 {a+1,b_{1},...,b_{k-1},b_k+1,b_{k+1},...,b_{n};} \hfill \\
 {c_{1} ,...,c_{k-1},c_k+1,c_{k+1},...,c_{n};} \hfill \\
\end{array}} y}\right]}
$
\begin{equation}
\label{eq25000}=F_{A}^{\left( {n} \right)} {\left[
{{\begin{array}{*{20}c}
 {a,b_{1},...,b_{n};} \hfill \\
 {c_{1} ,...,c_{k-1},c_k-1,c_{k+1},...,c_{n};} \hfill \\
\end{array}} y}\right]}-F_{A}^{\left( {n} \right)} {\left[ {{\begin{array}{*{20}c}
 {a,b_{1},...,b_{n};} \hfill \\
 {c_{1} ,...,c_{n};} \hfill \\
\end{array}} y}\right]}.
\end{equation}

Relations (\ref{eq24001})--(\ref{eq25000}) can be proved in two
ways: by comparing coefficients of equal powers of $y_1,$ ...,
$y_n$ on both sides or  mathematical induction.

For a given multiple hypergeometric function, it is useful to fund
a decomposition formula which would express the multivariable
hypergeometric function in terms of products of several simpler
hypergeometric functions involving fewer variables. Burchnall and
Chaundy \cite{{BC1},{BC2}} systematically presented a number of
expansions and decomposition formulas for some double
hypergeometric functions in series of simpler hypergeometric
functions. Using the Burchnall-Chaundy method, Hasanov and
Srivastava \cite{{HS6}, {HS7}} found decomposition formulas for a
whole class of hypergeometric functions in several (three and
more) variables. For example, the  Lauricella's hypergeometric
function $F_{A}^{\left( {n} \right)} $ defined by (\ref{eq7}) has
the decomposition formula \cite{HS6}
\begin{equation*}
F_{A}^{(n)} {\left[ {{\begin{array}{*{20}c}
 {a,b_{1} ,...,b_{n} ;} \hfill \\
 {c_{1} ,...,c_{n} ;} \hfill \\
\end{array}} y}  \right]}= {\sum\limits_{m_{2} ,...,m_{n} = 0}^{\infty}  {}} {\frac{{\left(
{a} \right)_{m_{2} + ... + m_{n}}  \left( {b_{1}}  \right)_{m_{2}
+ ... + m_{n}}  \left( {b_{2}}  \right)_{m_{2}}  ...\left( {b_{n}}
\right)_{m_{n}} }}{{m_{2} !...m_{n} !\left( {c_{1}} \right)_{m_{2}
+ ... + m_{n}}  \left( {c_{2}}  \right)_{m_{2}} ...\left( {c_{n}}
\right)_{m_{n}} }}}
\end{equation*}
\begin{equation*}
\times y_{1}^{m_{2} + ... + m_{n}} y_{2}^{m_{2}}
...y_{n}^{m_{n}}F{\left[ {{\begin{array}{*{20}c}
 {a + m_{2} + ... + m_{n} ,b_{1} + m_{2} + ... + m_{n};} \hfill \\
 {c_{1}+ m_{2} + ... + m_{n};} \hfill \\
\end{array}}y_{1}} \right]}
\end{equation*}
\begin{equation}\label{eqq11}
\times F_{A}^{(n-1)} {\left[ {{\begin{array}{*{20}c}
 {a + m_{2} + ... + m_{n} ,b_{2} + m_{2}
,...,b_{n} + m_{n};} \hfill \\
 {c_{2} + m_{2} ,...,c_{n} + m_{n};} \hfill \\
\end{array}} y_{2} ,...,y_{n}}  \right]},\,n \in {\mathbb{N}}\backslash {\left\{ {1}
\right\}}.
\end{equation}

However, due to the recurrence of formula (\ref{eqq11}),
additional difficulties may arise in the applications of this
expansion. Further study of the properties of the Lauricella's
hypergeometric functions showed that the formula (\ref{eqq11}) can
be reduced to a more convenient form.

\begin{lemma} \label{L11}
 \cite{A28}. The following decomposition formula
holds true
\begin{equation*}
 F_{A}^{(n)} {\left[
{{\begin{array}{*{20}c}
 {a,b_{1} ,...,b_{n} ;} \hfill \\
 {c_{1} ,...,c_{n} ;} \hfill \\
\end{array}} y}  \right]} = {\sum\limits_{{\mathop {m_{i,j} = 0}\limits_{(2 \le i \le j \le n)}
}}^{\infty}  {{\frac{{(a)_{A(n,n)}}} {M!}}}}\prod\limits_{k =
1}^{n}  \frac{{(b_{k} )_{B(k,n)}}}{(c_{k})_{B(k,n)}}
\end{equation*}
\begin{equation}\label{decomp}
\times
 \prod\limits_{k=1}^ny_{k}^{B(k,n)}F{\left[ {{\begin{array}{*{20}c}
 {a + A(k,n),b_{k} + B(k,n);} \hfill \\
 {c_{k} + B(k,n);} \hfill \\
\end{array}} y_{k}} \right]},\,\,n \in {\mathbb{N}}\backslash {\left\{ {1}
\right\}} ,\end{equation} \noindent where
\begin{equation*}
 M!:=m_{2,2}\cdot...\cdot m_{i,j} \cdot...\cdot
m_{n,n}, \,\,\,2 \le i \le j \le n;
\end{equation*}

\begin{equation*}
\label{eq13} A(k,n) = {\sum\limits_{i = 2}^{k + 1}
{{\sum\limits_{j = i}^{n} {m_{i,j}} }}} , \,\, B(k,n) =
{\sum\limits_{i = 2}^{k} {m_{i,k} +}}  {\sum\limits_{i = k +
1}^{n} {m_{k + 1,i}}}, \,\,k\in K.
\end{equation*}
\end{lemma}

\begin{lemma} \label{L12}
 If $a,b_{1},$ ..., $b_{n} $ are complex numbers
with $a \neq 0,\, - 1,\, - 2,...$ and $Re\left(a - b_{1} - ...-
 b_{n}\right)>0$,  then the following summation formula holds true
$$
{\sum\limits_{{\mathop {m_{i,j} = 0}\limits_{(2 \le i \le j \le
n)} }}^{\infty}{{\displaystyle\frac{{(a)_{A (n,n)}}}{{M!}}}}}
{\prod\limits_{k = 1}^{n} {{\frac{{\left( {b_{k}}\right)_{B(k,n)}
}}{{\left({a} \right)_{A(k,n)}}}\Gamma \left[{a - b_{k}}+{A(k,n) -
B(k,n)}\right]}}}
$$
\begin{equation}
\label{eq88899} ={\Gamma \left( {a - b_1-...-b_n}
\right)}{\Gamma^{n-1}(a)},\,\,n \in {\mathbb{N}}\backslash
{\left\{ {1} \right\}}.
\end{equation}
\end{lemma}

The Lemma \ref{L11} was proved by the method of mathematical
induction in \cite{A28}, the Lemma \ref{L12} is also proved
similarly.

From the formulae (\ref{decomp}) and (\ref{eq88899}) immediately
implies the following

\begin{corollary} \label{coro1}  If
$Re (a)> Re\left(b_1+...+b_n\right)>0$ and $Re \left(c_k\right)
>Re \left(b_k\right) >0$ \,($k\in K$), then the following
limiting equality
\begin{equation*}
{\mathop{\lim} \limits_{\mathop{y_{k} \to 0,}\limits_{k
=1,...,n}}}\left\{ y_{1}^{ - b_{1}}  ...y_{n}^{ - b_{n}}
F_{A}^{(n)} \left[\begin{array}{*{20}c}
 {a,b_{1} ,...,b_{n} ;} \hfill \\
 {c_{1} ,...,c_{n} ;} \hfill \\
\end{array}
1 -\frac{{1}}{{y_{1}}},...,1 - \frac{1}{y_n}\right]\right\}
\end{equation*}
\begin{equation}
\label{eq12222} = {\frac{{\Gamma \left( {a - b_1-...-b_n}
\right)}}{{\Gamma (a)}}}{\prod\limits_{k = 1}^{n} {{\frac{{\Gamma
\left( {c_{k}} \right)}}{{\Gamma \left( {c_{k} - b_{k}}
\right)}}}}}.
 \end{equation}
is valid.
\end{corollary}

One of the fundamental solutions of the equation (\ref{eq1}) that
we will use in this paper has the form \cite{A28}:
\begin{equation}
\label{eq2}q_{n} \left( {\xi,x}  \right) = \kappa_{n}
r^{-2\bar{\alpha}_n}(\xi x)^{(1-2\alpha)}
 F_{A}^{(n)} {\left[
{{\begin{array}{*{20}c}
 {\bar{\alpha}_n, 1 - \alpha_1,...,1-\alpha_n;} \hfill \\
 {2
- 2\alpha_1,...,2-2\alpha_n;} \hfill \\
\end{array}}\sigma}  \right]}
,
\end{equation}
where
$$
x:=\left(x_1,...,x_m\right), \,\,
\xi:=\left(\xi_1,...,\xi_m\right);\,\,(\xi
x)^{(1-2\alpha)}=\prod\limits_{k=1}^n\left(\xi_i
x_i\right)^{1-2\alpha_k};
$$
\begin{equation} \label{coeff}
\bar{\alpha}_n=\frac{m-2}{2}+\sum\limits_{k=1}^n\left(1-\alpha_k\right);\,\,\kappa_{n}
= 2^{2\bar{\alpha}_n-m}\frac{\Gamma \left( {\bar{\alpha}_n}
\right)}{{\pi ^{m /2}}}\prod\limits_{k=1}^n\frac{{\Gamma \left( {1
- \alpha_k} \right)}}{\Gamma \left( {2 - 2\alpha_k} \right)},
\end{equation}
\begin{equation}
\label{eqsigma} 0<2\alpha_k<1;
\,\,\sigma:=\left(\sigma_1,...,\sigma_n\right),\,\,\sigma_k = 1 -
{\frac{{r_{k}^{2}}} {{r^{2}}}},
\end{equation}
\begin{equation}
\label{eq4}  \quad r^{2} = {\sum\limits_{i = 1}^{m} {\left(
{\xi_{i} - x_{i}}  \right)^{2}}} , \quad r_{k}^{2} = \left(
{\xi_{k} + x_{k}}  \right)^{2} + {\sum\limits_{i = 1, i\neq k}^{m}
{\left( {\xi_{i} - x_{i}}  \right)^{2}}},\,\,k \in K .
\end{equation}

The fundamental solution given by (\ref{eq2}) possesses the
following potentially useful property:
\begin{equation}
\label{eq5} {\left. {q_{n} \left( \xi,x  \right)}
\right|}_{\xi_{k} = 0} = {\left. {q_{n} \left( {\xi,x}  \right)}
\right|}_{x_{k} = 0} = 0, \,\,k \in K .
\end{equation}

Throughout this paper, it is assumed that the dimension of the
space $m > 2.$

\section{Green's formula}

We consider the following identity:
\begin{equation}
\label{eqq7}  x^{(2\alpha)} \left[ uH_\alpha^{(m,n)} (v) -
vH_\alpha ^{(m,n)} (u) \right] = \sum\limits_{i = 1}^{m}
\frac{\partial} {\partial x_i}{\left[ x^{(2\alpha)} \left(
v_{x_{i}} u - vu_{x_i} \right)\right]} .
\end{equation}
Hereinafter,
\begin{equation}
\label{eq7000}
 x^{(2\alpha)}:= \prod\limits_{k=1}^nx_k^{2\alpha_k}.
 \end{equation}

 Integrating both sides of the identity (\ref{eqq7})
in a domain $\Omega$ located and bounded in $R_m^{n+}$, and using
the Gauss-Ostrogradsky formula, we obtain
\begin{equation}
\label{eq8000} \int_\Omega  x^{(2\alpha)} \left[ uH_\alpha
^{(m,n)} (v) - vH_\alpha^{(m,n)} (u) \right]dx  = \int_S
{{\left({uB_{Nx}^{\alpha}[v]-vB_{Nx}^{\alpha}[u]}\right)}} dS ,
\end{equation}
where $S\,$ is the boundary of $\Omega $, $N$ is the outer normal
to the surface $S$ and
\begin{equation}\label{e91}
B_{Nx}^{\alpha}[\,\,\,]=x^{(2\alpha)}\sum\limits_{i=1}^m\frac{\partial
[\,\,]}{\partial x_i} cos(N,x_i)
\end{equation}
is the conormal derivative with respect to $x$.

If $u$ and $v$ are solutions of the equation (\ref{eq1}), then we
find from the formula (\ref{eq8000}) that
\begin{equation}
\label{eq9} {\int_{S}
{{\left({uB_{Nx}^{\alpha}[v]-vB_{Nx}^{\alpha}[u]}\right)}}dS}  =
0.
\end{equation}

Assuming that $v = 1$ in (\ref{eq8000}) and replacing $u$ by
$u^{2}$, we obtain
\begin{equation}
\label{eq10} {\int_{\Omega}{x^{(2\alpha)}\sum\limits_{i = 1}^{m}
\left( {{\frac{{\partial u}}{{\partial x_{i}}} }} \right) ^{2}dx}}
= {\int_{S} u{B_{Nx}^{\alpha}[u]}dS} ,
\end{equation}
where $u\left({x}\right)$ is the solution of equation (\ref{eq1})

The special case of (\ref{eq9}) when $v = 1$ reduces to the
following form:
\begin{equation}
\label{eq11} {\int_{S} {B_{Nx}^{\alpha}[u]}dS}  = 0.
\end{equation}

We note from (\ref{eq11}) that the integral of the conormal
derivative of the solution of the equation (\ref{eq1}) along the
boundary $S$ of the domain $\Omega$ is equal to zero.

\section{A double-layer potential}

Let $\Gamma$ be a surface lying in $R_m^{n+}$ the boundary of
which on the hyperplane $x_k=0$ is denoted by $\gamma_k$ and
$\Omega $ be a finite domain in $R_{m}^{n+}  $, bounded by the
surface $\Gamma $  and the hyperplanes $x_{1} = 0$, ..., $x_n=0$.
The boundary of the domain $\Omega$ on the hyperplane $x_k=0$ is
denoted by $\Gamma_k$\,\,$(k \in K)$.

The surface $\Gamma$ in the Euclidean space $\,E_{m} $, that
satisfies the following three conditions is called the
\textit{Lyapunov surface} \cite{Mikh}:

(i).There is a definite normal at any point of the surface $\Gamma
$.

(ii).Let $x$ and $\xi $ be points of the surface $\Gamma $, and
$\vartheta $ angle between these normals. There exist positive
constants $a$ and $\kappa,$ such that
\begin{equation*}
\vartheta \leq a r^{\kappa} .
\end{equation*}

(iii).With respect to the surface $\Gamma $ we shall assume that
it approaches the hyperplanes $x_{1} = 0$, ..., $x_n=0$ under
right angle.

We introduce the following notation:
\begin{equation*} \label{not1}
\tilde{\xi}_k=\xi\setminus\{\xi_k\}=
\left(\xi_1,...,\xi_{k-1},\xi_{k+1},...,\xi_m\right),
\end{equation*}
\begin{equation*}\label{not2}
\xi^{(2\alpha)}=\prod\limits_{i=1}^n\xi_i^{2\alpha_i},
\,\,\,\tilde{\xi}_k^{(2\alpha)}=\prod\limits_{i=1,i\neq k}^n
\xi_i^{2\alpha_i},\,\,\,\tilde{\xi}_k^{(1)}=\prod\limits_{i=1,i\neq
k}^n \xi_i,
\end{equation*}
\begin{equation*}\label{not3}
x^{(1-2\alpha)}=\prod\limits_{i=1}^nx_i^{1-2\alpha_i};\,\,\,{X}_k^2=\xi_k^2+\sum\limits_{i=1,i\neq
k}^m\left(x_i-\xi_i\right)^2,\,\, k\in K.
\end{equation*}

We consider the following integral
\begin{equation}
\label{eq12} w^{(n)}\left(x\right) = {\int_{\Gamma} {\mu_{n}
\left( {\xi} \right)B_{N\xi}^{\alpha}[ q_{n} \left( {\xi ;x}
\right)]d_{\xi} \Gamma}}  ,
\end{equation}
where the density $\mu _{n} \left( {x} \right) \in C\left(
{\overline {\Gamma}}   \right)$, $q_{n} \left( {\xi ;x} \right)$
is given in (\ref{eq2}),  $N$ is the outer normal to the surface
$\Gamma$ and $B_{N\xi}^\alpha[\,\,\,]$ is the conormal derivative
with respect to $\xi$, as defined in (\ref{e91}).

We call the integral (\ref{eq12}) \textit{a double-layer potential
with density} $\mu _{n} \left( {\xi} \right)$.

We now investigate some properties of the double-layer potential
$w^{(n)} \left( {x} \right)$ when   $\mu _{n} \left( {\xi} \right)
\equiv 1$.

\begin{lemma} \label{L1} {The following formula holds true$:$}
\begin{equation*}
w_{1}^{(n)}(x) \equiv {\int_{\Gamma}{B_{N\xi}^{\alpha}[ q_{n}
\left( {\xi; x} \right)]d_{\xi} \Gamma}} = \left\{
{{\begin{array}{*{20}c}
 {i(x) - 1,\,\,\,\,x \in \Omega,} \hfill \\
 {i(x) - {\displaystyle\frac{{1}}{{2}}},\,\,\,\,x \in \Gamma ,} \hfill \\
 {\,\,i(x),\,\,\,\,\,\,\,\,\,x \notin {\Omega\cup\Gamma} ,} \hfill \\
\end{array}}}  \right.
\end{equation*}
where
\begin{equation*}
\begin{aligned}
i(x)&\equiv\sum\limits_{k=1}^n\int_{\Gamma_k}\tilde{\xi}_k^{(2\alpha)}\left.
\left(\xi_k^{2\alpha_k}\frac{\partial q_n(x;\xi)}{\partial
\xi_k}\right)\right|_{\xi_k=0}d_{\tilde{\xi}_k}\Gamma_k\\
&=\kappa_n
x^{(1-2\alpha)}\sum\limits_{k=1}^n(1-2\alpha_k)\int_{\Gamma_k}\tilde{\xi}_k^{(1)}X_{k}^{-2\bar{\alpha}_n}
\end{aligned}
\end{equation*}
\begin{equation} \label{ix}
\times F_{A}^{\left( {n-1} \right)} {\left[
{{\begin{array}{*{20}c}
 {\bar{\alpha} _{n} ,1 - \alpha _{1} ,...,1 - \alpha _{k-1}, 1 - \alpha_{k+1},...,1 - \alpha _{n} ;} \hfill \\
 {2 - 2\alpha _{1} ,...,2 - 2\alpha _{k-1}, 2 - 2\alpha_{k+1},2 - 2\alpha _{n} ;} \hfill \\
\end{array}} \tilde{\sigma}_{k}}\right]}d_{\tilde{\xi}_k}\Gamma_k,
\end{equation}

\begin{equation*} \label{not4}
\tilde{\sigma}_{k}=\left(-\frac{4x_1\xi_1}{{X}_k^2},...,-\frac{4x_{k-1}\xi_{k-1}}{{X}_k^2},-\frac{4x_{k+1}\xi_{k+1}}{{X}_k^2},...,-\frac{4x_n\xi_n}{{X}_k^2}\right).
\end{equation*}

{Here the domains $\Omega ,$ $\Gamma_k, $ $\gamma_k$ and the
surface $\Gamma $ are described as in this section.}
\end{lemma}

\textbf{Proof.}
We consider a few cases.\\

\textbf{Case 1.} When $x \notin {\Omega\cup\Gamma}$ and
$x_1>0,...,x_n>0$, it is noted that the function $q_{n} \left(
{\xi; x} \right)$ is a regular solution of the equation
(\ref{eq1}) inside of $\Omega$. Hence, in view of formula
(\ref{eq11}), we have
\begin{equation*}
w_{1}^{( {n})} \left( {x} \right) = i(x),\,x \notin
{\Omega\cup\Gamma},\,\,x_1>0,...,x_n>0,
\end{equation*}
where $i(x)$ is a function defined in (\ref{ix}).

\textbf{Case 2.} When $x \in \Omega ,$ we cut a ball centered at
$x$ with a radius $\rho $ off the domain $\Omega $ and denote the
remaining part by $\Omega _{\rho}  $ and the sphere of the
cut-off-ball by $C_{\rho}  $. The function $q_{n} \left( {\xi ; x}
\right)$ in (\ref{eq2}) is a regular solution of the equation
(\ref{eq1}) in the domain $\Omega_{\rho}$ and, in view of
(\ref{eq11}), we get following formula
\begin{equation}\label{eq16} w_{1}^{(n)} (x) =i(x)+ {\mathop {\lim} \limits_{\rho \to
0}} {\int_{C_{\rho}} {B_{N\xi}^\alpha{\left[{q_{n} \left( {\xi; x}
\right)} \right]}d_{\xi} C_{\rho}} },
\end{equation}

We calculate  $B_{N\xi}^\alpha{\left[{q_{n} \left( {\xi; x}
\right)} \right]}$ with respect to $\xi$. Using the formula of
differentiation (\ref{eq24001}), adjacent relations (\ref{eq25}),
(\ref{eq25000}) and a definition of the conormal derivative
(\ref{e91}), we obtain (for details, see \cite{Turk})
\begin{equation}
\label{eq2611}B_{N\xi}^\alpha{\left[{q_{n} \left( {\xi; x}
\right)} \right]} =B_1(\xi; x)B_{N\xi}^\alpha{\left[ {\ln
\frac{1}{r}}\right]}+B_2(\xi; x),
\end{equation}
where
\begin{equation}
\label{eq26} B_1(\xi; x)= 2\bar{\alpha}_{n}
\kappa_{n}\frac{\left(\xi x\right)^{(1-2\alpha)}}{r^{2\bar
{\alpha}_{n}}} F_{A}^{\left( {n} \right)} {\left[
{{\begin{array}{*{20}c}
 {1+\bar{\alpha} _{n},1 - \alpha _{1} ,...,1 - \alpha _{n} ;} \hfill \\
 {\,\,\,2 - 2\alpha _{1}\,\, ,...,\,\,2 - 2\alpha _{n} ;} \hfill \\
\end{array}} \sigma}  \right]},
\end{equation}

$B_2(\xi; x)= \kappa_{n} x^{(1-2\alpha)}
 r^{-2\bar{\alpha}_{n}}\sum\limits_{k = 1}^{n}{{{\left(1 -
2\alpha_{k}\right)}} {{\tilde{\xi}_{k}^{(1)}}}}\cos \left( {N,
\xi_{k}} \right) $
\begin{equation}
\label{eq2600}  \times F_{A}^{\left( {n} \right)} {\left[
{{\begin{array}{*{20}c}
 {\bar{\alpha} _{n},1 - \alpha _{1} ,...,1 - \alpha _{n} ;} \hfill \\
 {2 - 2\alpha _{1} ,...,2 - 2\alpha _{k-1}, 1-2\alpha_{k},2-2\alpha_{k+1},...,2-2\alpha_n ;} \hfill \\
\end{array}} \sigma}  \right]}.
\end{equation}

It is easy to see  that
\begin{equation}
\label{eqnull}{\mathop {\lim} \limits_{\xi_k \to
0}}B_{N\xi}^\alpha{\left[{q_{n} \left( {\xi; x} \right)} \right]}
={\mathop {\lim} \limits_{x_k \to 0}}B_{N\xi}^\alpha{\left[{q_{n}
\left( {\xi; x} \right)} \right]}=0, \,\, k\in K.
\end{equation}

Taking (\ref{eq2611}) into account we first calculate the
following integral
\begin{equation*}
 j_{1}(x;\rho)= \int_{C_{\rho}} B_{1}
(\xi;x)B_{N\xi}^\alpha{\left[{\ln \frac{1}{r}}\right]}dC_{\rho}.
\end{equation*}

We use the following generalization spherical system of
coordinates:
\begin{equation}  \label{eq18}
\xi_i=x_i+\rho\Phi_i(\varphi), \,\,\, i=1,...,m, \,\,
\end{equation}
where
\begin{equation*}
\varphi=\left(\varphi_1,...,\varphi_{m-1}\right);
\,\,\Phi_1(\varphi) = \cos\varphi_{1},
\end{equation*}
\begin{equation*}
\Phi_i(\varphi)=\sin \varphi _{1} \sin \varphi _{2} \cdots\sin
\varphi _{i - 1} \cos \varphi _{i} ,i = 2,...,m - 1,
\end{equation*}
\begin{equation*}
\Phi_m(\varphi)=\sin \varphi _{1} \sin \varphi _{2} \cdots\sin
\varphi _{m - 2} \sin \varphi _{m - 1}
\end{equation*}
\begin{equation*}
{(0 \leq \rho \leq r,\,\,\,0 \leq \varphi _{1} \leq \pi
,\,\,...,\,\, 0 \leq \varphi _{m - 2} \leq \pi ,\,\,\,0 \leq
\varphi _{m - 1} \leq 2\pi )}.
\end{equation*}

Then we have
\begin{equation*}
 j_{1} (x,\rho)  = -2\bar{\alpha}_n \kappa_{n}{\int\limits_{0}^{2\pi} {d\varphi _{m -
1}}}\int\limits_{0}^{\pi}  \sin \varphi _{m - 2} d\varphi _{m -
2}...{\int\limits_{0}^{\pi}}{\sin^{m - 2}\varphi
_{1}}\aleph(\rho,\varphi)d\varphi _{1},
\end{equation*}
where
\begin{equation*}
  \aleph(\rho,\varphi)=  \prod\limits_{k = 1}^{n}\frac{x_k}{r_{k\rho}}\cdot \prod\limits_{k = 1}^{n}
  \left(\frac{x_k+\rho \Phi_k(\rho)}{r_{k\rho}}\right)^{1-2\alpha_k}\cdot \prod\limits_{k = 1}^{n} \left(\frac{\rho^2}{r_{k\rho}^2}\right)^{1-\alpha_k}
\end{equation*}
\begin{equation}\label{sumalef}
 \times F_{A}^{\left( {n} \right)} {\left[
{{\begin{array}{*{20}c}
 {1+\bar{\alpha} _{n},1 - \alpha _{1} ,...,1 - \alpha _{n} ;} \hfill \\
 {\,\,\,2 - 2\alpha _{1}\,\, ,...,\,\,2 - 2\alpha _{n} ;} \hfill \\
\end{array}} 1-\frac{\rho^2}{r_{1\rho}^2},...,1-\frac{\rho^2}{r_{n\rho}^2}}
\right]},
\end{equation}
where
\begin{equation*}
r_{k\rho}^2 =4x_k^2+4\rho x_k\Phi_k(\varphi) +\rho^2, \,\,k\in K.
\end{equation*}

It is easy to see that when $\rho \to 0$ the function
$\aleph(\rho,\varphi) $ becomes an expression that does not depend
on $\varphi$. Indeed,  applying the formula (\ref{eq12222}), we
get
\begin{equation} \label{alef}{\mathop{\lim} \limits_{\rho \to 0}}\aleph(\rho,\varphi)
=2^{-2\bar{\alpha}_n+m-2}\frac{\Gamma \left({{{{m}}/{{2}}}}
\right)}{\Gamma\left(1+\bar{\alpha}_n\right)}\prod\limits_{k=1}^n\frac{\Gamma\left(2-2\alpha_k\right)}{\Gamma\left(1-\alpha_k\right)}.
\end{equation}

It is not difficult to establish that
\begin{equation*}
{\int\limits_{0}^{2\pi}  {d\varphi _{m - 1}}}
{\int\limits_{0}^{\pi} {\sin \varphi _{m - 2} d\varphi _{m - 2}}}
{\int\limits_{0}^{\pi}  {\sin ^{2}\varphi _{m - 3} d\varphi _{m -
3}}}  ...
\end{equation*}
\begin{equation} \label{eq33333}
  ...{\int\limits_{0}^{\pi}  {\sin ^{m - 2}\varphi _{1}
d\varphi _{1}}} =\frac{2\pi^{m/2}}{\Gamma\left(m/2\right)},
\,\,m=2,3,....
\end{equation}

If we take into account (\ref{sumalef}), (\ref{alef}),
(\ref{eq33333}) and (\ref{coeff}), then we will have
\begin{equation} \label{minus1}
 {\mathop{\lim} \limits_{\rho \to 0}}j_1(x,\rho)=-1.
\end{equation}

By similar evaluations one can get that
\begin{equation} \label{nuli}
 {\mathop{\lim} \limits_{\rho \to 0}}j_2(x,\rho)=0.
\end{equation}

Substituting (\ref{minus1}) and (\ref{nuli}) into (\ref{eq16}), we
finally get
\begin{equation*}
 w_{1}^{(n)}(x) = i(x)-1.
\end{equation*}

\textbf{Case 3.} When $x \in \Gamma ,$ we cut a sphere $C_{\rho} $
centered at $x$ with a radius $\rho $ off the domain $\Omega $ and
denote the remaining part of the surface by ${\Gamma}_\rho ',$
that is, ${\Gamma}_\rho ' = \Gamma - \Gamma _{\rho} .$ Let
${C}'_{\rho}  $ denote a part of the sphere $C_{\rho}  $ lying
inside the domain $\Omega .$ We consider the domain $\Omega
_{\rho}  $ which is bounded by a surface ${\Gamma}_\rho '$,\, a
semisphere ${C}'_{\rho} $, hyperplanes $\Gamma _{1}
$,...,$\Gamma_n$ and their boundaries $\gamma_1 $,..., $\gamma_n$.
Then we have
\begin{equation}
\label{eq23} w_{1}^{\left( {n} \right)} \left( {x} \right) =i(x)+
{\mathop {\lim }\limits_{\rho \to 0}} \int_{{\Gamma}_\rho '}
B_{N\xi}^\alpha\left[ q_{n} \left( {\xi; x} \right)\right]d_{\xi}
{\Gamma}_\rho' .
\end{equation}

When the point $x$ lies outside the domain $\Omega _{\rho}  $, it
is found that, in this domain $q_{n} \left( {\xi; x} \right)$ is a
regular solution of the equation (\ref{eq1}). Therefore, by virtue
of (\ref{eq11}), we have
\begin{equation}
\label{eq24} {\int_{{\Gamma}_\rho'} {B_{N\xi}^\alpha\left[q_{n}
\left({\xi ; x} \right)\right]d_{\xi} {\Gamma}_\rho'}}  =
{\int_{{C}'_{\rho}}{B_{N\xi}^\alpha\left[q_{n}\left({\xi; x}
\right)\right]d_{\xi} {C}'_{\rho}}}.
\end{equation}

Substituting (\ref{eq24}) into (\ref{eq23}), we get
\begin{equation*}
w_{1}^{\left( {n} \right)} \left( {x} \right) = i(x)+ {\mathop
{\lim }\limits_{\rho \to 0}} {\int_{C_{\rho} ^{'}}
{B_{N\xi}^\alpha\left[q_{n}\left({\xi;x} \right)\right]d_{\xi}
{C}'_{\rho}}}. \end{equation*}

Similarly, by again introducing the spherical coordinates
(\ref{eq18}) centered at the point $x$, we have
\begin{equation*}
w_{1}^{({n})} \left( {x} \right) = i(x) - {\frac{{1}}{{2}}}, x \in
\Gamma .
\end{equation*}

\textbf{Case 4.} Finally, we put the point $x$ on the any
hyperplane $x_k=0\,(k\in K)$.  By virtue of (\ref{eqnull}), the
equality $\left. w_{1}^{\left( {n} \right)}(x)\right|_{x_k=0}=0\,
(k \in K)$ immediately follows from the definition (\ref{eq12}).

The proof of Lemma \ref{L1} is  completed.

\begin{theorem} \label{T1} If $x \in \Gamma$,   then the following inequality
holds true$:$
\begin{equation}
\label{eq27} {\left| B_{N\xi}^\alpha\left[{q_{n} \left( {\xi; x}
\right)}\right] \right|} \leq {\frac{C_1}{{r^{m -
2}}{r^{(2-2\alpha)}}}},
\end{equation}
where  $m > 2$ and $\alpha:=\left(\alpha_1,...,\alpha_n\right)$
are real parameters with $0 < 2\alpha_k < 1$ $\,(k \in K)$ as in
the equation $(\ref{eq1});$  $r$ and $r^{(2-2\alpha)} $ are as in
(\ref{eq4}) and (\ref{eq1213}) respectively, \,\,$C_1$ is a
constant.
\end{theorem}

\textbf{Proof.} We transform a right side of the equality
(\ref{eq26}). Making use of the decomposition formula
(\ref{decomp}), we have
\begin{gather}\label{eq12333}
\begin{aligned}
  B_1(\xi; x)&= 2\bar{\alpha}_{n} \kappa_{n}\frac{\left(\xi
x\right)^{(1-2\alpha)}}{r^{2\bar {\alpha}_{n}}}\\
&\times{{\sum\limits_{{\mathop {m_{i,j} = 0}\limits_{(2 \le i \le
j \le n)} }}^{\infty}}
{{\frac{{\left(1+\bar{\alpha}_n\right)_{A(n,n)}}}
{M!}}}}\prod\limits_{k = 1}^{n}  \frac{{\left(1-\alpha_{k}
\right)_{B(k,n)}}}{\left(2-2\alpha_{k}\right)_{B(k,n)}}\\
&\times
 \prod\limits_{k=1}^n \sigma_{k}^{B(k,n)}F{\left[ {{\begin{array}{*{20}c}
 {1+\bar{\alpha}_n + A(k,n),1-\alpha_{k} + B(k,n);} \hfill \\
 {2-2\alpha_{k} + B(k,n);} \hfill \\
\end{array}} \sigma_{k}} \right]}.
\end{aligned}
\end{gather}

 Applying the formula (\ref{trans}) to
each Gaussian hypergeometric function in (\ref{eq12333}), we
obtain
\begin{gather}\label{eq1211}
\begin{aligned}
&B_1(\xi; x) = \frac{2\bar{\alpha}_{n} \kappa_{n}\left(\xi
x\right)^{(1-2\alpha)}}{r^{m-2}r^{(2-2{\alpha})}}\\
 &\times
{\sum\limits_{{\mathop {m_{i,j} = 0}\limits_{(2 \le i \le j \le
n)} }}^{\infty} {{\frac{{\left(1+\bar{\alpha}_n\right)_{A(n,n)}}}
{M!}}}}\prod\limits_{k = 1}^{n}  \frac{{\left(1-\alpha_{k}
\right)_{B(k,n)}}}{\left(2-2\alpha_{k}\right)_{B(k,n)}}\left(-\omega_k\right)^{B(k,n)}\\
&\times
 \prod\limits_{k=1}^n F{\left[ {{\begin{array}{*{20}c}
 {1-2\alpha_k-\bar{\alpha}_n +B(k,n)-A(k,n),1-\alpha_{k} + B(k,n);} \hfill \\
 {2-2\alpha_{k} + B(k,n);} \hfill \\
\end{array}} \omega_k}  \right]},
\end{aligned}
\end{gather}
where
\begin{equation}\label{eq1213}
r^{(2-2\alpha)}:=\prod\limits_{k=1}^n
r_k^{2-2\alpha_k},\,\,\,\,{\omega}_k:=1-\displaystyle\frac{r^2}{r_k^2},\,\,k\in
K.
\end{equation}

Similarly, from the formula (\ref{eq2600}) we find
\begin{gather}\label{eq1212}
\begin{aligned}
  &B_2(\xi; x)= \frac{\kappa_{n}x^{(1-2\alpha)}}{r^{m-2}r^{(2-2{\alpha})}}\sum\limits_{k=1}^n \tilde{\xi}_k^{(1)}\cos\left(N,\tilde{\xi}_k\right)
   {\sum\limits_{{\mathop {m_{i,j} = 0}\limits_{(2 \le i
\le j \le n)} }}^{\infty}
{{\frac{\left(\bar{\alpha}_n\right)_{A(n,n)}} {M!}}}}\\
 &\times
\left(1-2\alpha_k+B(k,n)\right)
 \prod\limits_{i = 1}^{n}  \frac{{\left(1-\alpha_{i}
\right)_{B(i,n)}}}{\left(2-2\alpha_{i}\right)_{B(i,n)}}\left(-\omega_i\right)^{B(i,n)}\\
&\times
 \prod\limits_{i=1,i\neq k}^n F{\left[ {{\begin{array}{*{20}c}
 {2-2\alpha_i-\bar{\alpha}_n +B(i,n)-A(i,n),1-\alpha_{i} + B(i,n);} \hfill \\
 {2-2\alpha_{i} + B(i,n);} \hfill \\
\end{array}} \omega_i}
\right]}\\
&\times F{\left[ {{\begin{array}{*{20}c}
 {1-2\alpha_k-\bar{\alpha}_n +B(k,n)-A(k,n),1-\alpha_{k} + B(k,n);} \hfill \\
 {1-2\alpha_{k} + B(k,n);} \hfill \\
\end{array}} \omega_k}
\right]}.
\end{aligned}
\end{gather}

Proof of the inequality (\ref{eq27}) follows from the formulae
(\ref{eq2611}),(\ref{eq1212}) and (\ref{eq1211}).

\begin{theorem}\label{T2} {If a surface $\Gamma $
satisfies the conditions $(i)$--$(iii),$ then the following
inequality holds true$:$
\begin{equation} \label{eq1215}
{\int_{\Gamma}  {{\left|B_{N\xi}^\alpha\left[ q_{n} \left( {\xi ;
x} \right)\right] \right|}d_{\xi} \Gamma}} \leq C_2,
\end{equation}
where $C_2$  is a constant.}
\end{theorem}

\textbf{Proof.} The inequality (\ref{eq1215}) immediately follows
from the formulae (\ref{eq2611}), (\ref{eq1212}) and
(\ref{eq1211}).

\begin{theorem} \label{T3}
The following limiting formulas hold true for a double-layer
potential $(\ref{eq12})$$:$
\begin{equation}
\label{eq28} w_{i}^{\left( {n} \right)} \left( {t} \right) = -
{\frac{{1}}{{2}}}\mu _{n} \left( {t} \right) + {\int_{\Gamma} {\mu
_{n} \left( {s} \right)K_{n} \left( {s;t} \right)d_{s} \Gamma}}
\end{equation}
and
\begin{equation}
\label{eq29} w_{e}^{\left( {n} \right)} \left( {t} \right) =
{\frac{{1}}{{2}}}\mu _{n} \left( {t} \right) + {\int_{\Gamma} {\mu
_{n} \left( {s} \right)K_{n} \left( {s;t} \right)d_{s} \Gamma}}  ,
\end{equation}
where \begin{equation*}\mu _{n} \left( {t} \right) \in C\left(
{\overline {\Gamma}} \right), s: = \left( {s_{1} ,s_{2}
,\cdots,s_{m}} \right),t: = \left( {t_{1} ,t_{2} ,\cdots,t_{m}}
\right),
\end{equation*}
\begin{equation*}
 K_{n} \left( {s;t} \right) = B_{Ns}^\alpha\left[
q_{n} \left( {s; t} \right)\right],\,\, s \in \Gamma,\,\,  t\in
\Gamma,
\end{equation*}
$w_{i}^{\left( {n} \right)} \left( {t} \right)$  and
$w_{e}^{\left( {n} \right)} \left( {t} \right)$ are limiting
values of the double-layer potential (\ref{eq12}) at the point $t
\in \Gamma $ from the inside and the outside$,$ respectively.
\end{theorem}

\textbf{Proof.} Theorem \ref{T3} follows from Lemma \ref{L1} and
Theorem \ref{T2}.

\section{A simple-layer potential}

We consider the following integral
\begin{equation}
\label{eq31} v_n\left( {x} \right) = {\int_{\Gamma}  {\rho _{n}
\left( {\xi}  \right)q_{n} \left( {\xi ; x} \right)d_{\xi}
\Gamma}}  ,
\end{equation}
where the density $\rho _{n} \left( {x} \right) \in C
\left(\overline {\Gamma}\right)  $ and $q_{n} \left( {\xi; x}
\right)$ is given in (\ref{eq2}). We call the integral
(\ref{eq31}) \textit{a simple-layer potential with density} $\rho
_{n} \left( {\xi} \right)$.

The simple-layer potential (\ref{eq31}) is defined throughout the
domain  $R_m^{n+}$ and a continuous function when passing through
the surface $\Gamma.$ Obviously, a simple-layer potential
$v_n\left( {x} \right)$ is a regular solution of the equation
(\ref{eq1}) in any domain lying in the  $R_m^{n+}$. It is easy to
see that as the point $x$ tends to infinity, a simple-layer
potential $v_n(x)$ tends to zero. Indeed, we let the point $x$ be
on the $2^{n}$th part of the sphere given by
\begin{equation*}
{\sum\limits_{i = 1}^{m} {x_{i}^{2}}}= R^{2},x_{1} > 0,...,x_n>0.
\end{equation*}

Now we consider fundamental solution $q_n(\xi,x)$ defined in
(\ref{eq2}). Applying consequentially the decomposition formula
(\ref{decomp}) and formula (\ref{trans}), we obtain
\begin{gather}\label{eq121111}
\begin{aligned}
  &q_n(\xi; x)= \frac{\kappa_{n}\left(\xi
x\right)^{(1-2\alpha)}}{r^{m-2}r^{(2-2{\alpha})}}\\
&\times{\sum\limits_{{\mathop {m_{i,j} = 0}\limits_{(2 \le i \le j
\le n)} }}^{\infty}
{{\frac{{\left(\bar{\alpha}_n\right)_{A(n,n)}}}
{M!}}}}\prod\limits_{k = 1}^{n}  \frac{{\left(1-\alpha_{k}
\right)_{B(k,n)}}}{\left(2-2\alpha_{k}\right)_{B(k,n)}}\left(-\omega_k\right)^{B(k,n)}\\
&\times
 \prod\limits_{k=1}^n F{\left[ {{\begin{array}{*{20}c}
 {2-2\alpha_k-\bar{\alpha}_n +B(k,n)-A(k,n),1-\alpha_{k} + B(k,n);} \hfill \\
 {2-2\alpha_{k} + B(k,n);} \hfill \\
\end{array}} \omega_k}  \right]}.
\end{aligned}
\end{gather}

Then, by virtue of (\ref{eq121111}), we have
\begin{equation*}
{\left| {v_n(x)} \right|} \leq {\int_{\Gamma}  {{\left| {\rho _{n}
(t)} \right|}{\left| {q_{n} (t;x)} \right|}d_{t} \Gamma}} \leq
M_1R^{2-m-n}\,\,\,(R \geq R_{0} ),
\end{equation*}
where $ M_1 $ is a constant.

We take an arbitrary point $P(x_{0})$ on the surface $\Gamma $ and
draw a normal $N$ at this point. Consider on this normal any point
$M(x)$, not lying on the surface $\Gamma $, we find the conormal
derivative of the simple-layer potential (\ref{eq31}):
\begin{equation}
\label{e301} B_{Nx}^\alpha\left[ v_n(x)\right] = {\int_{\Gamma}
{B_{Nx}^\alpha\left[ q_{n} \left( {\xi; x} \right)\right]d_{\xi}
\Gamma}}.
\end{equation}

The integral (\ref{e301}) exists also in the case when the point
$M(x)$ coincides with the point $P,$ which we mentioned above.

\begin{theorem} \label{T4}
The following limiting formulas hold true for a simple-layer
potential $(\ref{eq31})$$:$
\begin{equation}
\label{eq32}B_{Nt}^\alpha\left[ v_n(t)\right]_i =
{\frac{{1}}{{2}}}\rho _{n} \left( {t} \right) + {\int_{\Gamma}
{\rho _{n} \left( {s} \right)K_{n} \left( {s;t} \right)d_{s}
\Gamma}}
\end{equation}
and
\begin{equation}
\label{eq33} B_{Nt}^\alpha\left[ v_n(t)\right]_e = -
{\frac{{1}}{{2}}}\rho _{n} \left( {t} \right) + {\int_{\Gamma}
{\rho _{n} \left( {s} \right)K_{n} \left( {s;t} \right)d_{s}
\Gamma}}  ,
\end{equation}
where
\begin{equation*} \rho _{n}  \in C\left(
{\overline {\Gamma}}   \right), \,\,K_{n} \left( {s;t} \right) =
B_{Nt}^\alpha\left[ q_n(s;t)\right],\,\, s\in \Gamma,\,\,
t\in\Gamma,
\end{equation*}
$B_{Nt}^\alpha\left[ v_n(t)\right]_i $ and $B_{Nt}^\alpha\left[
v_n(t)\right]_e$ are limiting values of the normal derivative of
simple-layer potential (\ref{eq31}) at the point $t \in \Gamma $
from the inside and the outside, respectively.
\end{theorem}

Making use of these formulas, the jump on the normal derivative of
the simple-layer potential follows immediately:
\begin{equation}
\label{eq34} B_{Nt}^\alpha\left[
v_n(t)\right]_i-B_{Nt}^\alpha\left[ v_n(t)\right]_e =
 \rho _{n} (t).
\end{equation}

By virtue of formulae (\ref{eq1211}) and  (\ref{eq1212}), one can
prove that, when the point $x$ tends to infinity, the following
inequality:
\begin{equation*}
{\left| B_{Nt}^\alpha\left[ v_n(t)\right]\right|} \leq M_2R^{ -
2\bar{\alpha}_n}\,\,\,(R \geq R_{0} )
\end{equation*}
is valid, where $M_2$ is a constant.

In exactly the same way as in the derivation of (\ref{eq10}), it
is not difficult to show that Green's formulas are applicable to
the simple-layer potential (\ref{eq31}) as follows:
\begin{equation}
\label{eq35} {\int_{\Omega}  {x^{(2\alpha)}  {\sum\limits_{k =
1}^{m} {\left( {{\frac{{\partial v_n}}{{\partial x_{k}}} }}
\right)}} dx}}  = {\int_{S} v_n(x){B_{Nx}^\alpha\left[
v_n(x)\right]_i dS}}
\end{equation}
and
\begin{equation}
\label{eq36} {\int_{\Omega'}  {x^{(2\alpha)}  {\sum\limits_{k =
1}^{m} {\left( {{\frac{{\partial v_n}}{{\partial x_{k}}} }}
\right)}} dx}}  = - {\int_{S}v_n(x) {B_{Nx}^\alpha\left[
v_n(x)\right]_e dS}} ,
\end{equation}
Hereinafter  $\Omega':=R_{m}^{n+} \setminus \overline {\Omega}$
\,\, is a infinite domain.  $x^{(2\alpha)}$  is defined by the
formula (\ref{eq7000}).

\section{Integral equations for densities}\label{S5}

Formulas (\ref{eq28}), (\ref{eq29}), (\ref{eq32}) and (\ref{eq33})
can be written as the following integral equations for densities:
\begin{equation}
\label{eq37} \mu _{n} \left( {s} \right) - \lambda {\int_{\Gamma}
{K_{n} \left( {s;t} \right)\mu _{n} \left( {t} \right)d_{t}
\Gamma}}   = f_{n} (s)
\end{equation}
and
\begin{equation}
\label{eq38} \rho _{n} \left( {s} \right) - \lambda {\int_{\Gamma}
{K_{n} \left( {t;s} \right)\rho _{n} \left( {t} \right)d_{t}
\Gamma}}   = g_{n} (s),
\end{equation}
where
\begin{equation*}
\lambda = 2,\,\, f_{n} (s) = - 2w_{i}^{\left( {n} \right)} \left(
{s} \right), \,\,g_{n} (s) = - 2B_{Ns}^\alpha \left[ {v_n(s)}
\right]_e
\end{equation*}
and
\begin{equation*}
\lambda = - 2,\,\, f_{n} (s) = 2w_{e}^{\left( {n} \right)} \left(
{s} \right),\,\, g_{n} (s) = 2B_{Ns}^\alpha
\left[{v_n(s)}\right]_i.
\end{equation*}

Equations (\ref{eq37}) and (\ref{eq38}) are mutually conjugated
and, by Theorem \ref{T1}, Fredholm theory is applicable to them.
We show that $\lambda = 2$ is not an eigenvalue of the kernel
$K_{n} \left( {s,t} \right)$. This assertion is equivalent to the
fact that the homogeneous integral equation
\begin{equation}
\label{eq39} \rho _{n} \left( {t} \right) - 2{\int_{\Gamma} {K_{n}
\left( {s;t} \right)\rho _{n} \left( {s} \right)d_{s} \Gamma}} =0
\end{equation}
has no non-trivial solutions.

Let $\rho_{0} \left( {t} \right)$ be a continuous non-trivial
solution of the equation (\ref{eq39}). The simple-layer potential
with density $\rho_{0} \left( {t} \right)$ gives us a function
$v_{0} \left( {x} \right),$ which is a solution of the equation
(\ref{eq1}) in the domains $\Omega $ and $\Omega'$. By virtue of
the equation (\ref{eq39}), the limiting values of the conormal
derivative of $B_{Ns}^\alpha\left[v_{0} (s)\right]_{e}$ are zero.
The formula (\ref{eq36}) is applicable to the simple-layer
potential $v_{0} (x)$, from which it follows that $v_{0} (x) =
const$ in domain ${\Omega'} $. At infinity, a simple layer
potential is zero, and consequently $v_{0} (x) \equiv 0$ in
${\Omega'} $, and also on the surface $\Gamma.$ Applying now
(\ref{eq35}), we find that $v_{0} (x) \equiv 0$ is valid also
inside the region $\Omega.$ But then $B_{Ns}^\alpha\left[v_{0}
(s)\right]_{i} = 0, $ and by virtue of formula (\ref{eq34}) we
obtain $\rho_{0} \left( {t} \right) \equiv 0.$ Thus, clearly, the
homogeneous equation (\ref{eq39}) has only the trivial solution;
consequently, $\lambda = 2$ is not an eigenvalue of the kernel
$K_{n} \left( {s;t} \right)$.

Similarly, one can show that $\lambda = 2$ is not an eigenvalue of
the kernel $K_{n} \left( {s,t} \right)$.

\section{The uniqueness of the solution of Dirichlet's problem}

We apply the obtained results of potential theory to the solving
the boundary value problem for the equation (\ref{eq1}) in the
domain $\Omega$.

We introduce the following notation:
\begin{equation*} \label{not4}
\tilde{x}_k=x\setminus\{x_k\}=
\left(x_1,...,x_{k-1},x_{k+1},...,x_m\right),\,\,\tilde{x}_k^{(2\alpha)}=\prod\limits_{i=1,i\neq
k}^n x_i^{2\alpha_i},\,\,k\in K.
\end{equation*}

\textbf{The Dirichlet problem.} Find a regular solution of the
equation (\ref{eq1}) in the domain $\Omega$  that is continuous in
the closed domain $\overline {\Omega}$ and satisfies the following
boundary conditions:
\begin{equation}
\label{e401}
 {\left. {u} \right|}_{\Gamma}  = \varphi (x),\,\,\,x \in \overline{\Gamma},
\end{equation}
\begin{equation} \label{e400} \left.u(x) \right|_{x_k=0}= \tau_k\left(\tilde{x}_k\right),\,\,\,\tilde{x}_k \in
\overline{\Gamma}_{k},\, k \in K,
\end{equation}
where $\varphi (x)$ and $\tau_k\left(\tilde{x}_k\right)$ are given
continuous functions fulfilling the following matching conditions:
$ \left. \varphi(x)\right|_{\Gamma_k}=\left.
\tau_k\left(\tilde{x}_k\right)\right|_\Gamma, \,\, k \in K.$

 Considering the equality
(\ref{eq10}), we obtain
\begin{gather}\label{eq41}
\begin{aligned}
 {\int_{\Omega}  {x^{(2\alpha)}  {\sum\limits_{i =
1}^{m} {\left( {{\frac{{\partial u}}{{\partial x_{i}}} }}
\right)^2}}dx}}&=- \int_{\Gamma}\varphi(x)B_{Nx}^\alpha[u]
d_x\Gamma\\
& + \sum\limits_{k=1}^n{\int_{\Gamma_{k}}
{\tilde{x}_k^{(2\alpha)}\tau_k\left(\tilde{x}_k\right)
\left.\left(x_{k}^{2\alpha}  {\frac{{\partial u}}{{\partial
x_{k}}}}\right)\right|_{x_k=0} d_{\tilde{x}_k}\Gamma_k}} .
\end{aligned}
\end{gather}

In case of the homogeneous Dirichlet  problem from (\ref{eq41})
one can easily get
\begin{equation*}
{\int_{\Omega}  {x^{(2\alpha)}  {\sum\limits_{i = 1}^{m} {\left(
{{\frac{{\partial u}}{{\partial x_{i}}} }} \right)}} ^{2}dx}}  =
0.
\end{equation*}

Hence, it follows that $u(x) = 0 $ in $\overline {\Omega}  $.

Thus we have proved the following

\begin{theorem} If the Dirichlet problem has a regular solution,  then it is
unique.
\end{theorem}

\section{Green's function revisited}

To solve this problem, we use the Green's function method. First,
we construct the Green's function for solving the Dirichlet
problem for an equation (\ref{eq1}) in the domain $\Omega$ bounded
by an arbitrary surface $\Gamma$ and hyperplanes $x_1=0,$ ...,
$x_n=0$. In the end, we show that, thanks to the Green's function,
the solution of the Dirichlet problem in a special domain (in the
$2^n$th part of the multidimensional ball) takes a simpler form.

\noindent \textbf{Definition.} We refer to $G_n(x;\xi )_{} $ as
Green's function
of the Dirichlet problem, if it satisfies the following conditions:\\

\noindent {\bf Condition 1.} The function $G_n(x;\xi )_{}$ is a
regular solution of the equation (\ref{eq1}) in the domain
$\Omega,$ except at the point $\xi ,$ which is any fixed
point of $\Omega$.\\

\noindent {\bf Condition 2.} The function $G_n(x;\xi )_{}$
satisfies the boundary conditions given by
\begin{equation}
\label{eq42} {\left. {G_n(x;\xi )} \right|}_{\Gamma}  = 0, {\left.
{G_n(x;\xi )} \right|}_{x_{k} = 0} = 0, \,\,k \in K.
\end{equation}

\noindent {\bf Condition 3.} The function $G_n(x;\xi )_{}$ can be
represented as follows:
\begin{equation}
\label{eq43} G_n(x;\xi ) = q_{n} (x;\xi ) + v_{n} (x;\xi ),
\end{equation}
where  $q_{n} \left( {x;\xi}  \right)$ is a fundamental solution
of the equation (\ref{eq1}), defined in (\ref{eq2}) and the
function $v_{n} (x;\xi )$ is a regular solution of the equation
(\ref{eq1}) in the domain $\Omega$.

The construction of the Green's function $G_n(x,\xi )_{}$ reduces
to finding its regular part $v_{n} (x;\xi )$ which, by virtue of
(\ref{eq42}), (\ref{eq43}) and (\ref{eq5}), must satisfy the
following boundary conditions:
\begin{equation}
\label{eq44} {\left. {v_{n} (x;\xi )} \right|}_{\Gamma}  = -
{\left. {q_{n} (x;\xi )} \right|}_{\Gamma}
\end{equation}
and
\begin{equation*}
{\left. {v_{n} (x;\xi )} \right|}_{x_{k} = 0} = 0, \,\,k \in K.
\end{equation*}

We look for the function $v_{n}(x;\xi )$ in the form of a
double-layer potential given by
\begin{equation}
\label{eq45} v_{n} \left( {x;\xi}  \right) = {\int_{\Gamma}
{\mu_{n} \left( {t;\xi}\right)B_{Nt}^\alpha \left[q_{n} \left(
{t;x} \right)\right]d_{t}\Gamma}} .
\end{equation}

Taking into account the equality (\ref{eq28}) and the boundary
condition (\ref{eq44}), we obtain the integral equation for the
density  $\mu _{n} \left( {s;\xi}  \right)$ as follows:
\begin{equation}
\label{eq46} \mu _{n} \left( {s;\xi}  \right) - 2{\int_{\Gamma}
{\mu _{n} \left( {t;\xi } \right)K_{n} \left( {s;t} \right)d_{t}
\Gamma}}   = 2q_{n} (s;\xi ),\,\,\,\,s \in \Gamma .
\end{equation}

The right-hand side of (\ref{eq46}) is a continuous function with
respect to $s$ (the point $\xi $ lies inside $\Omega )$.  By
Theorem \ref{T1}, Fredholm theory is applicable to the equation
(\ref{eq46}). In section \ref{S5} it was proved that $\lambda=2$
is not an eigenvalue of the kernel $K_n(s,t)$ and, consequently,
the equation (\ref{eq46}) is solvable and its continuous solution
can be written in the following form:
\begin{equation}
\label{eq47} \mu _{n} \left( {s;\xi}  \right) = 2q_{n} (s;\xi ) +
4{\int_{\Gamma}  {R_{n} \left( {s,t;2} \right)q_{n} \left( {t;\xi}
\right)d_{t} \Gamma}}  ,
\end{equation}
where $R_{n} (s,t;2)$ is the resolvent of the kernel $K_{n}
(s;t)$, $s \in \Gamma $. Substituting (\ref{eq47}) into
(\ref{eq45}), we obtain
\begin{equation*}
 {v_{n} \left( {x,\xi}  \right) = 2{\int_{\Gamma}
{q_{n} (t;\xi )B_{Nt}^\alpha \left[q_{n} \left( {t;x}
\right)\right]d_{t} \Gamma}}  }
\end{equation*}
\begin{equation}
\label{eq48}
 + 4{\int_{\Gamma} {{\int_{\Gamma}  { R_{n} \left(
{t,s;2} \right)q_{n} \left( {s;\xi}  \right)B_{Nt}^\alpha
\left[q_{n} \left( {t;x} \right)\right]d_{t} \Gamma}}  d_{s}
\Gamma }} .
\end{equation}

We now define a following function:
\begin{equation}
\label{eq49} g\left( {x} \right) = {\left\{
{{\begin{array}{*{20}c}
 {v_{n} (x;\xi ),\,\,x \in \Omega ,} \hfill \\
 { - q_{n} (x;\xi ),\,\,x \in {\Omega' }.} \hfill \\
\end{array}}}  \right.}
\end{equation}

The function $g\left( {x} \right)$ is a regular solution of
(\ref{eq1}) both inside the domain $\Omega ,$ and inside ${\Omega'
}$ and equal to zero at infinity. Since point $\xi $ lies inside
$\Omega,$ then in $ {\Omega' }$ the function $g\left( {x} \right)$
has derivatives of any order in all variables, continuous up to
$\Gamma .$ We can consider $g\left( {x} \right)$ in ${\Omega' }$
as a solution of the equation (\ref{eq1}) satisfying the boundary
conditions given by
\begin{equation}
\label{eq50} {\left.B_{Nx}^\alpha \left[g\left({x}\right)\right]
\right|}_{\Gamma}  = - {\left. B_{Nx}^\alpha \left[q_{n} \left(
{x;\xi} \right)\right]\right|}_{\Gamma}
\end{equation}
and
\begin{equation*}
{\left.g\left({x} \right)\right|}_{x_{k}=0}=0, \,\,k \in K.
\end{equation*}

We represent this solution in the form of a simple-layer potential
as follows:
\begin{equation}
\label{eq51} g\left( {x} \right) = {\int_{\Gamma}  {\rho _{n}
\left( {t;\xi} \right)q_{n} \left( {t;x} \right)d_{t} \Gamma
,\,\,\,x \in  {\Omega' }}}
\end{equation}
with an unknown density $\rho_{n} \left( {t;\xi}  \right)$.

Using the formula (\ref{eq33}), by virtue of condition
(\ref{eq50}), we obtain the following integral equation for the
density $\rho _{n} \left( {t;\xi}  \right)$:
\begin{equation}
\label{eq52} \rho _{n} \left( {s;\xi}  \right) - 2{\int_{\Gamma}
{\rho _{n} \left( {t;\xi}  \right)K_{n} \left( {t;s} \right)d_{t}
\Gamma}}   = 2B_{Ns}^\alpha \left[q_{n} \left( {s;\xi}
\right)\right].
\end{equation}
Equation (\ref{eq52}) is conjugated with the equation
(\ref{eq46}). Its right-hand side is a continuous function with
respect to $s$. Thus, clearly, the equation (\ref{eq52}) has the
following continuous solution:
\begin{equation}
\label{eq53} \rho _{n} \left( {s;\xi}  \right) = 2B_{Ns}^\alpha
\left[q_{n} \left( {s;\xi} \right)\right]  + 4{\int_{\Gamma}
{R_{n} \left( {t,s;2} \right)B_{Nt}^\alpha \left[q_{n} \left(
{t;\xi} \right)\right]d_{t}\Gamma}}.
\end{equation}

According to (\ref{eqnull})  is obvious that the function $\rho
_{n} \left( {s;\xi} \right)$  has the following limiting values
\begin{equation}
\label{eqlimit}{\mathop {\lim} \limits_{\xi_k\to 0}} \rho _{n}
\left( {s;\xi}  \right) = 0, \,\,s \in \Gamma, \,\, k \in K.
\end{equation}

 The values of a simple-layer potential $g(x)_{}$ on the
surface $\Gamma $ are equal to $ - q_{n} \left( {x;\xi}  \right)$,
that is, just as the values of the function $v_{n} \left( {x;\xi}
\right)$ and  on the hyperplanes $x_{k} = 0$, $k\in K$   are equal
to zero. Hence, by virtue of the uniqueness theorem for the
Dirichlet problem, it follows that the formula (\ref{eq51}) for
the function $g(x)$ defined by (\ref{eq49}) holds throughout in
the  $x_{1} \geq 0,$...,$x_{n} \geq 0,$  that is,
\begin{equation}
\label{eq54} v_{n} \left( {x;\xi}  \right) = {\int_{\Gamma}  {\rho
_{n} \left( {t;\xi} \right)q_{n} \left( {t;x} \right)d_{t} \Gamma
,\,\,\,x \in \Omega .}}
\end{equation}

Thus, the regular part $v_{n} (x;\xi )$ of Green's function is
representable in the form of a simple-layer potential.

Applying the formula (\ref{eq32}) to (\ref{eq54}), we obtain
\begin{equation*}
2B_{Ns}^\alpha \left[v_{n} \left( {s;\xi} \right)\right]_i = \rho
_{n} \left( {s;\xi}  \right) + 2{\int_{\Gamma}  {K_{n} \left(
{t;s} \right)\rho _{n} \left( {t;\xi}  \right)d_{t} \Gamma}}  ,
\end{equation*} but, according to (\ref{eq52}), we have
\begin{equation*}
2B_{Ns}^\alpha \left[q_{n} \left( {s;\xi} \right)\right]_i = \rho
_{n} \left( {s;\xi}  \right) - 2{\int_{\Gamma}  {K_{n} \left(
{t;s} \right)\rho _{n} \left( {t;\xi}  \right)d_{t} \Gamma}}  .
\end{equation*}

Summing the last two equalities by term-wise and taking into
account (\ref{eq43}), we have
\begin{equation}
\label{eq55} B_{Ns}^\alpha \left[G_{n} \left( {s;\xi}
\right)\right] = \rho _{n} \left( {s;\xi}  \right),
\end{equation}
and, consequently, formula (\ref{eq54}) can be written in the
following form:
\begin{equation*}
 v_{n} \left( {x;\xi}  \right) = {\int_{\Gamma} {q_{n}
\left( {t;x} \right) B_{Nt}^\alpha \left[G_{n} \left( {t;\xi}
\right)\right]d_{t} \Gamma .}}
\end{equation*}

Multiplying both sides of (\ref{eq53}) by $q_{n} \left( {s;x}
\right)$, integrating by  $s$ over the surface $\Gamma $ and, by
virtue of (\ref{eq47}) and (\ref{eq45}), we obtain

$$
v_{n} \left( {\xi ;x} \right) = {\int_{\Gamma}  {\rho _{n} \left(
{t;\xi} \right)q_{n} \left( {t;x} \right)d_{t} \Gamma .}}
$$
Comparing this last equation with the formula (\ref{eq54}), we
have
\begin{equation}
\label{eq57} v_{n} \left( {\xi ;x} \right) = v_{n} \left( {x;\xi}
\right),
\end{equation}
if the points $x$ and $\xi$ are inside the domain $\Omega.$

\begin{lemma} \label{L3}
If points $x $  and $\xi $ are inside domain $\Omega,$ then
Green's function $G_{n} \left( {x;\xi} \right)$  is symmetric
about those points.
\end{lemma}

\textbf{Proof.} The proof of Lemma \ref{L3} follows from the
representation (\ref{eq43}) of Green's function and the equality
(\ref{eq57}).

For a  domain $\Omega_{0}$ bounded by the hyperplanes $x_{1} =
0$,..., $x_n=0$ and the  $2^n$th part of the sphere given by
$$
x_{1}^{2} + x_{2}^{2} + \cdots + x_{m}^{2} = R^{2},\,\,
x_1>0,...,x_n>0,
$$
Green's function of the Dirichlet problem has the following form:
\begin{equation}
\label{eq58} G_{0n} \left( {x;\xi}  \right) = q_{n} \left( {x;\xi}
\right) - \left( {{\frac{{R}}{{\varrho}}}}
\right)^{2\bar{\alpha}_n} \cdot q_{n} \left(
{x;{\frac{{R^{2}}}{{\varrho^{2}}}}\xi} \right),
\end{equation}
where
\begin{equation}
\label{eqroro}
 \varrho^{2}= \xi _{1}^{2} + \xi _{2}^{2} + \cdots + \xi _{m}^{2}.
\end{equation}

We show that the function given by
$$
v_{0n} \left( {x;\xi}  \right) = - \left(
{{\frac{{R}}{{\varrho}}}} \right)^{2\bar{\alpha}_n}\cdot q_{n}
\left( {x;{\frac{{R^{2}}}{{\varrho^{2}}}}\xi}\right)
$$
can be represented in the following form:
\begin{equation}
\label{v02} v_{0n} \left( {x;\xi} \right) = - {\int_{\Gamma} {\rho
_{n} \left( {s;x} \right)v_{0n} \left( {s;\xi}  \right)d_{s}
\Gamma ,}}
\end{equation}
where $\rho _{n} \left( {s;x} \right)$ is a solution of the
equation (\ref{eq52}).

Indeed, by letting an arbitrary point $\xi $ be inside the domain
$\Omega $, we consider the function given by
\begin{equation*}
\label{eq59} u\left( {x;\xi}  \right) = - {\int_{\Gamma}  {\rho
_{n} \left( {s;x} \right)v_{0n} \left( {s;\xi}  \right)d_{s}
\Gamma .}}
\end{equation*}

The function  $u(x,\xi)$ satisfies the equation (\ref{eq1}), since
this equation is satisfied by the function $\rho _{n} (s;x)$.
Substituting the expression (\ref{eq53}) for $\rho _{n} (s;x)$, we
obtain
\begin{equation}
\label{eq60} u\left( {x;\xi}  \right) = - {\int_{\Gamma}  {\psi
\left( {s;\xi} \right)B_{Ns}^\alpha \left[q_{n} \left( {s;x}
\right)\right]d_{s} \Gamma ,}}
\end{equation}
where
$$
\psi \left( {s;\xi}  \right) = 2v_{0n} \left( {s;\xi}  \right) +
4{\int_{\Gamma}  {R_{n} \left( {s,t;2} \right)v_{0n} \left(
{t,\xi} \right)d_{t} \Gamma}}  ,
$$
that is, $\psi \left( {s;\xi}  \right)$ is a solution of the
integral equation
\begin{equation}
\label{eq61} \psi \left( {s;\xi}  \right) - 2{\int_{\Gamma} {K_{n}
\left( {s;t} \right)\psi \left( {t;\xi}  \right)d_{t} \Gamma}}   =
2v_{0n} \left( {s;\xi } \right).
\end{equation}

Applying the formula (\ref{eq28}) to the double-layer potential
(\ref{eq60}), we obtain
$$
u_{i} \left( {s;\xi}  \right) = {\frac{{1}}{{2}}}\psi \left(
{s;\xi} \right) - {\int_{\Gamma}  {K_{n} \left( {s;t} \right)\psi
\left( {t;\xi} \right)d_{t} \Gamma}}  ,
$$
whence, by virtue of (\ref{eq61}), we get
$$
u_{i} \left( {s;\xi}  \right) = v_{0n} \left( {s;\xi}  \right), s
\in \Gamma .
$$

It is easy to see that
$$
{\left. u(x;\xi ) \right|}_{x_{k} = 0} = 0, \,\,\,{\left.
v_{0n}(x;\xi)\right|}_{x_{k} = 0} = 0, \,\,k\in K.
$$
Thus, clearly, the functions $u\left( {x;\xi}  \right)$ and
$v_{0n} \left( {x;\xi} \right)$ satisfy the same equation
(\ref{eq1}) and the same boundary conditions. Also, by virtue of
the uniqueness of the solution of the Dirichlet problem, the
equality
$$u\left({x;\xi} \right) \equiv v_{0n} \left({x;\xi}\right)$$ is
satisfied.

Now, subtracting the expression (\ref{eq58}) from (\ref{eq43}), we
obtain
$$
H_{n} \left( {x;\xi} \right) = G_{n} \left({x;\xi}  \right) -
G_{0n} \left( {x;\xi} \right) = v_{n} \left({x;\xi} \right) -
v_{0n} \left( {x;\xi} \right)
$$
or, by virtue of (\ref{eq54}), (\ref{eq58}), (\ref{v02})
 and  a symmetry property (\ref{eq57}), we obtain
\begin{equation}
\label{eq62} H_{n} \left( {x;\xi}  \right) = {\int_{\Gamma}
{G_{0n} \left( {t;\xi} \right)\rho _{n} \left( {t;x} \right)d_{t}
\Gamma}} .
\end{equation}

\section{Solving the Dirichlet problem for equation
(\ref{eq1})}

Let $\xi$ be a point inside the domain $\Omega$. Consider the
domain $\Omega_{\varepsilon,\delta}\subset \Omega$, where
$\delta:=\left(\delta_1,...,\delta_n\right)$, bounded by the
surface $\Gamma_{\varepsilon}$ which is parallel to the surface
$\Gamma,$ and the domains $\Gamma_{k\delta_k}$ lying on the
hyperplanes $x_{k}=\delta_k>\varepsilon,\,k \in K.$ We choose
$\delta_1$,..., $\delta_n$ and  $\varepsilon $  so small that the
point $x_0$ is inside $\Omega _{\varepsilon,\delta}$. We cut out
from the domain $\Omega _{\varepsilon ,\delta}$ a ball of small
radius $\rho$ with center at the point $x_{0}$ and the remainder
part of $\Omega _{\varepsilon,\delta}$ denote by $\Omega
_{\varepsilon ,\delta} ^{\rho} ,$ in which the Green's function
$G_{n}\left({x;x_{0}} \right)$ is a regular solution of the
equation (\ref{eq1}).

Let $u(x)$ be a regular solution of the equation (\ref{eq1}) in
the domain $\Omega $ that satisfies the boundary conditions
(\ref{e401}) and (\ref{e400}). Applying the formula (\ref{eq9}),
we obtain
$$
{\int_{\Gamma _{\varepsilon}}   \left(
G_{n}B_{Nx}^\alpha\left[u\right] - u B_{Nx}^\alpha\left[G_n\right]
\right)d_x\Gamma _{\varepsilon}} +
\sum\limits_{k=1}^n{\int_{\Gamma _{k\delta_k}}\left(
uB_{Nx}^\alpha\left[G_n\right]\right.}
$$
$${\left.\left. - G_n B_{Nx}^\alpha\left[u\right]
\right)\right|_{x_k=\delta_k}d_{x_k}\Gamma _{k\delta_k}} =
{\int_{C_{\rho}}{\left( G_{n}B_{Nx}^\alpha\left[u\right] - u
B_{Nx}^\alpha\left[G_n\right] \right)d_xC_{\rho}}} .
$$

Passing to the limit as $\rho \to 0$ and then as $\delta_1 \to
0$,..., $\delta_n \to 0$ and $\varepsilon \to 0$, we obtain
\begin{equation*}
 u(\xi) = \sum\limits_{k=1}^n
{\int_{\Gamma_{k}}{\tau_k \left( {\tilde{x}_k}
\right)\left.\left(x_k^{2\alpha_k}\frac{\partial
G_{n}\left({x;\xi}\right)}{\partial
x_k}\right)\right|_{x_k=0}d_{\tilde{x}_k}\Gamma_k}}
\end{equation*}
\begin{equation}
\label{eq63} -{\int_{\Gamma} {\varphi (x){B_{Nx}^\alpha
\left[G_n\left(x;\xi\right)\right] }d_{x}
\Gamma}}=\sum\limits_{k=1}^nT_k\left(\xi\right)+\Phi\left(\xi\right).
\end{equation}

We show that the formula (\ref{eq63}) gives a solution to the
Dirichlet problem.

It is easy to see that each integral $T_{k}(\xi)$ in the formula
(\ref{eq63}) satisfies the  equation (\ref{eq1}) and is regular in
the domain $\Omega $, continuous in $\overline {\Omega}$.

We use the following notations:
\begin{equation}
\label{eq6464}
\begin{aligned}
\vartheta_k (\xi) &= {\int_{\Gamma _{k}} {\tau_k\left(
{\tilde{x}_k}\right)\tilde{q}_n(\tilde{x}_{k0};\xi)d_{\tilde{x}_k}\Gamma_k}}\\
&= {\int_{\Gamma _{k}} {\tilde{x}_k^{(2\alpha)}\tau_k\left(
{\tilde{x}_k}\right)\left.\left(x_k^{2\alpha_k}\frac{\partial
q_{n}\left({x;\xi}\right)}{\partial
x_k}\right)\right|_{x_k=0}d_{\tilde{x}_k}\Gamma_k}}
\end{aligned}
\end{equation}
where
$$
\tilde{x}_{k0}=\left(x_1,...,x_{k-1},0,x_{k+1},...,x_m\right),\,\,k\in
K.$$

Here $\vartheta_k(\xi)$ is a continuous function in $\overline
{\Omega}$. In view of (\ref{eq6464}) and (\ref{eq48}) and the
symmetry property of the function $v_{n}(x;\xi)$, the integral
$T_{k}(\xi)$ can be represented in the following form:
\begin{equation*}
 T_{k}(\xi)=\vartheta_k(\xi)+2{\int_{\Gamma}
{\vartheta_k(t)B_{Nt}^\alpha
\left[q_n\left(t;\xi\right)\right]d_{t}\Gamma}}
\end{equation*}
\begin{equation}
\label{eq65} + 4{\int_{\Gamma}{{\int_{\Gamma}{R_{n}
(t,s;2)\vartheta_k(s)B_{Nt}^\alpha
\left[q_n\left(t;\xi\right)\right]d_{t}\Gamma d_{s}\Gamma,}}}
}\,\,k\in K.
\end{equation}

The last two integrals in the formula (\ref{eq65}) are
double-layer potentials. Taking the formula (\ref{eq28}) and the
integral equation for the resolvent $R_{n} (t,s;2)$ into account,
we obtain
$$
{\left. {T_{k}(\xi)}\right|}_{\Gamma}=0,\,\,k\in K.
$$

It is easy to prove that
$$
{\mathop {\lim} \limits_{\xi_k\to 0}}T_{k}(\xi) = \tau
(\tilde{\xi}_k), \,\,\tilde{\xi}_k \in \overline{\Gamma}_{k}, \,\,
k\in K.
$$
Indeed, by virtue of (\ref{eq54}) and the symmetry property of the
function $v_n(x,\xi)$, the integral $T_k(\xi)$ can be rewritten as
\begin{equation*}
T_{k}(\xi)={\int_{\Gamma_k}{\tau_k\left( {\tilde{x}_k}
\right)\tilde{q}_n(\tilde{x}_{k0};\xi)d_{\tilde{x}_k}\Gamma_k}}
\end{equation*}
\begin{equation}
\label{eq6500}+{\int_{\Gamma_k}{\tau_k\left( {\tilde{x}_k}
\right)d_{\tilde{x}_k}\Gamma_k
\int_\Gamma\rho_n\left(t;\xi\right)\tilde{q}_n(\tilde{x}_{k0};t)d_t\Gamma}},
\end{equation}
where $\rho_n\left(t;\xi\right)$  is defined in (\ref{eq53}).

Following the work \cite{Sm}, one can get that
$$
{\mathop {\lim} \limits_{\xi_k\to 0}}{\int_{\Gamma_k}{\tau_k\left(
{\tilde{x}_k}
\right)\tilde{q}_n(\tilde{x}_{k0};\xi)d_{\tilde{x}_k}\Gamma_k}}
 = \tau_k\left( {\tilde{\xi}_k}
\right), \,\,{\tilde{\xi}_k} \in \overline{\Gamma}_{k},\,\, k\in K
.
$$
Taking (\ref{eqlimit}) into account we see that second addend in
(\ref{eq6500}) is  zero at $\xi_k = 0,$\, $k\in K$.

Now consider the last integral $\Phi(\xi)$ in the formula
(\ref{eq63}), which, by virtue of (\ref{eq53}) and (\ref{eq55}),
can be written in the following form:
$$
\Phi(\xi) = -{\int_{\Gamma}{\varphi(s)\rho_{n}(s;\xi)d_{s}\Gamma
}}  = - {\int_{\Gamma}  {\theta (t)B_{Nt}^\alpha
\left[q_n\left(t;\xi\right)\right]d_{t}\Gamma}},
$$
where
$$
\theta (t) = 2\varphi(t)+ 4{\int_{\Gamma}\varphi(s){R_{n}
(t,s;2)d_{s}\Gamma }} ,
$$
that is, the function $\theta(s)$ is a solution to the integral
equation
\begin{equation}
\label{eq66} \theta (s) - 2{\int_{\Gamma}  {K_{n} (s;t)\theta
(t)d_{t} \Gamma = 2\varphi (s).}}
\end{equation}

Since the function $\theta (s)$ is continuous, the function
$\Phi(\xi)$ is a solution to the equation (\ref{eq1}), regular in
the domain $\Omega $ and continuous in $\overline {\Omega}$ which,
by virtue of (\ref{eq28}), (\ref{eq29}) and (\ref{eq66}),
satisfies the following condition:
$$
{\left.{\Phi(\xi)}\right|}_{\Gamma}  = \varphi (s), \,\,s\in
\overline{\Gamma}.
$$

It is easy to see that
\begin{equation*}
{\mathop {\lim} \limits_{\xi_k\to 0}}\Phi(\xi) = 0,
\,\,\tilde{\xi}_k \in \overline{\Gamma}_{k}, \,\, k\in K.
\end{equation*}

The formula (\ref{eq63}), and with it all the proof, requires that
$m>2$. However, the formula (\ref{eq63}) is also valid for $m=2$
(in case of $m=2$, for details, see \cite{Sm}).

Thus we have proved the following
\begin{theorem} If $\tau_k\left(
{\tilde{x}_k} \right)\in C\left(\overline{\Gamma}_k\right)\cap
C^2\left({\Gamma}_k\right)$  and  $\varphi(x)\in
C\left(\overline{\Gamma}_k\right)\cap C^2\left({\Gamma}_k\right)$
are given functions fulfilling the matching conditions $ \left.
\varphi(x)\right|_{\Gamma_k}=\left.
\tau_k\left(\tilde{x}_k\right)\right|_\Gamma, \,\, k \in K$  ,
then the Dirichlet problem for equation $H_\alpha^{(m,n)}(u)=0$\,
$(m\geq 2, 0<n\leq m)$ in the domain $\Omega$ has unique solution
represented by formula (\ref{eq63}).
\end{theorem}

Using the formulae (\ref{eq58}) and (\ref{eq62}) we rewrite a
solution (\ref{eq63}) to the Dirichlet problem for equation
(\ref{eq1}) in the following form
\begin{equation*}
 u(\xi) = \sum\limits_{k=1}^n
{\int_{\Gamma_{k}}{\tilde{x}_k^{(2\alpha)}\tau_k \left(
{\tilde{x}_k}
\right)\left[\tilde{G}_{0n}\left(\tilde{x}_{k0};\xi\right)+\tilde{H}_{n}\left(\tilde{x}_{k0};\xi\right)\right]d_{\tilde{x}_k}\Gamma_k}}
\end{equation*}
\begin{equation}
\label{eq63121} -{\int_{\Gamma} {\varphi (x)\left\{{B_{Nx}^\alpha
\left[G_{0n}\left(x;\xi\right)\right]+B_{Nx}^\alpha
\left[H_n\left(x;\xi\right)\right] }\right\}d_{x} \Gamma}},
\end{equation}
where
\begin{equation*}
\tilde{G}_{0n}\left(\tilde{x}_{k0};\xi\right) =
(1-2\alpha_k)\kappa_n \xi^{(1-2\alpha)}
\tilde{x}_k^{(1-2\alpha)}\left[\frac{\tilde{F}_{A}^{(n-1)}\left(\tilde{\sigma}_{k}\right)}{{X}_k^{2\bar{\alpha}_n}}-\frac{\tilde{F}_{A}^{(n-1)}\left(\tilde{\omega}_{k}\right)}{{Y}_k^{2\bar{\alpha}_n}}\right],
\end{equation*}

\begin{equation*}
\tilde{x}_{k0}=\left(x_1,...,x_{k-1},0,x_{k+1},...,x_m\right);
\,\,{X}_k^2=\xi_k^2+\sum\limits_{i=1,i\neq
k}^m\left(x_i-\xi_i\right)^2,
\end{equation*}

\begin{equation*}
Y_k^2= {{\sum\limits_{i = 1, i\neq k}^{m} {\left( {R -
{\displaystyle\frac{{x_{i} \xi_{i}}} {{R}}}} \right)^{2} +
{\displaystyle\frac{{1}}{{R^{2}}}}}} {\sum\limits_{i = 1,i\neq
k}^{m} {x_{i}^{2} {\sum\limits_{j = 1,j \ne i}^{m} {\xi_{j}^{2}}}
} }  - (m - 2)R^{2}},
\end{equation*}

\begin{equation*}
\tilde{\sigma}_{k}=\left(-\frac{4x_1\xi_1}{{X}_k^2},...,-\frac{4x_{k-1}\xi_{k-1}}{{X}_k^2},-\frac{4x_{k+1}\xi_{k+1}}{{X}_k^2},...,-\frac{4x_n\xi_n}{{X}_k^2}\right),
\end{equation*}

\begin{equation*}
\tilde{\omega}_{k}=\left(-\frac{R^2}{\varrho^2}\frac{4x_1\xi_1}{{Y}_k^2},...,-\frac{R^2}{\varrho^2}\frac{4x_{k-1}\xi_{k-1}}{{Y}_k^2},-\frac{R^2}{\varrho^2}\frac{4x_{k+1}\xi_{k+1}}{{Y}_k^2},...,-\frac{R^2}{\varrho^2}\frac{4x_n\xi_n}{{Y}_k^2}\right),
\end{equation*}

\begin{equation*}
\tilde{F}_{A}^{(n-1)}\left(...\right)= F_{A}^{(n-1)} {\left[
{{\begin{array}{*{20}c}
 {\bar{\alpha}_n,1-\alpha_{1},...,1-\alpha_{k-1},1-\alpha_{k+1},...,1-\alpha_n;} \hfill \\
 {2-2\alpha_{1},...,2-2\alpha_{k-1},2-2\alpha_{k+1},...,2-2\alpha_n;} \hfill \\
\end{array}} ...}  \right]},
\end{equation*}

\begin{equation*}
\xi^{(1-2\alpha)}=\prod\limits_{i=1}^n\xi_i^{1-2\alpha_i},\,\,{x}^{(1-2\alpha)}=\prod\limits_{i=1}^n
x_i^{1-2\alpha_i},
\end{equation*}

\begin{equation*}
\tilde{H}_{n}\left(\tilde{x}_{k0};\xi\right)=\left.\left(x_k^{2\alpha_k}\frac{\partial
H_{n}(x,\xi)}{\partial x_k}\right)\right|_{x_k=0},\,\,k\in K;
\end{equation*}

\begin{equation*}
 H_{n} \left( {x;\xi}  \right) = {\int_{\Gamma}
{G_{0n} \left( {t;\xi} \right)\rho _{n} \left( {t;x} \right)d_{t}
\Gamma}};\,\,\,m\geq 2, 0<n\leq m.
\end{equation*}
Here $\bar{\alpha}_n$ and $\kappa_n$ are defined in (\ref{coeff}).

 We remark that the solution (\ref{eq63121}) to the Dirichlet
problem is more convenient for further investigations. The
resulting explicit integral representation (\ref{eq63121}) plays
an important role in the study of problems for equation of the
mixed type (that is, elliptic-hyperbolic or elliptic-parabolic
types): it makes it easy to derive the basic functional
relationship between the traces of the sought solution and of its
derivative on the line of degeneration from the elliptic part of
the mixed domain.

\section{A solution in case of the $2^n$th part of the $m$-dimensional ball}

In this section we find the solution to the Dirichlet problem for
equation (\ref{eq1}) in the special domain $\Omega_0$  bounded by
the $2^n$th part of the sphere:
$$
S=\left\{x: \sum\nolimits_{i=1}^m x_i^2 =R^2,\,\,\, x_j>0,\,\,j\in
K\right\}
$$
and
\begin{equation*} S_k=\left\{x: \, \sum\nolimits_{i=1,i\neq k}^m
x_i^2 \leq R^2,\,\,\, x_j>0,\,\, j\in
K\setminus\{k\}\right\},\,\,k\in K.
\end{equation*}

In case of the domain $\Omega_{0}$, the function $H_{n}(x;\xi)
\equiv 0$ and the solution (\ref{eq63121}) assumes a simpler form:
\begin{equation}\label{eq63129}
\begin{aligned}
 u(\xi) =& \kappa_n
\xi^{(1-2\alpha)}\sum\limits_{k=1}^n\left(1-2\alpha_k\right)\\
&\times{\int_{S_{k}}{\tilde{x}_k^{(1)}
\left[\frac{\tilde{F}_{A}^{(n-1)}\left(\tilde{\sigma}_{k}\right)}{{X}_k^{2\bar{\alpha}_n}}-\frac{\tilde{F}_{A}^{(n-1)}\left(\tilde{\omega}_{k}\right)}{{Y}_k^{2\bar{\alpha}_n}}\right]\tau_k
\left( {\tilde{x}_k} \right)d_{\tilde{x}_k}S_k}}\\
 &+2\bar{\alpha}_n\kappa_n\xi^{(1-2\alpha)}\\
 &\times{\int_{S}
{x^{(1)}F_{A}^{(n)} {\left[ {{\begin{array}{*{20}c}
 {1+\bar{\alpha}_n,1-\alpha_{1},...,1-\alpha_n;} \hfill \\
 {2-2\alpha_{1},...,2-2\alpha_n;} \hfill \\
\end{array}} \sigma}  \right]}\frac{R^2-\varrho^2}{Rr^{2+2\bar{\alpha}_n}}\varphi(x)
d_{x} S}},
\end{aligned}
\end{equation}
where
\begin{equation*}
{x}^{(1)}=\prod\limits_{i=1}^n
x_i,\,\,\,\tilde{x}_k^{(1)}=\prod\limits_{i=1,i\neq k}^n
x_i,\,\,k\in K.
\end{equation*}
Here $\sigma$, $r$ and $\varrho$ are defined in (\ref{eqsigma}),
(\ref{eq4}) and (\ref{eqroro}) respectively.

The formula (\ref{eq63129}) was found by  other way in
\cite{Turk}, but, of course, here we are interested in obtaining
this formula as an application (example) of the potential theory
constructed in the present paper.

We note that particular  cases of the formula (\ref{eq63129}) for
the two- and three-dimensional singular elliptic equations were
known \cite{{A17},{A18},{SI},{Sm}}.

\section{\bf Concluding remarks and observations}

In such  widely-investigated subject as potential theory, both
simple-layer  potential and double-layer potential play
significant role in solving boundary value problems involving
various families of elliptic partial differential equations.  In
particular, a double-layer potential provides a solution of
Laplace's  equation corresponding to the electrostatic or magnetic
potential associated with a dipole distribution on a closed
surface in the $m$-dimensional Euclidean space.

In our present investigation of the multidimensional singular
elliptic equation (\ref{eq1}), we use potential theory results in
order to represent boundary value problems in integral equation
form. In fact, in problems with known Green's functions, an
integral equation formulation leads to powerful numerical
approximation schemes. Thus, by seeking the representation of the
solution of the boundary value problem as a double-layer potential
with unknown density, we are eventually led to a Fredholm equation
of the second kind for the explicit determination of the solution
in terms of hypergeometric functions in many variables.
Lauricella's  hypergeometric function $F_A^{(n)}(a, b_1,...,b_n;
c_1,...,c_n; y_1,...,y_n)$ possesses easily-accessible numerical
algorithms for computational purposes, can indeed be used to
numerically compute the solution presented here for many different
special values of the parameters $a$, $b_k,$ $c_k$ and of the
arguments $z_k$, $\,\,k\in K$.

Numerical applications of several suitably specialized versions of
the solutions presented in this paper can be found in solid
mechanics, fluid mechanics, elastic dynamics, electro-magnetics,
and acoustics (see, for details, some of the citations
\cite{{Bers},{Frankl2}} handling special situations which were
motivated by such widespread applications).
\\

\end{document}